\documentclass{article}
\usepackage[a4paper, total={6in, 8in}]{geometry}
\usepackage[utf8]{inputenc}
\usepackage{amsmath, amsfonts,amsthm}
\usepackage{float}
\usepackage{cite}
\usepackage{graphicx}
\usepackage{textcomp}
\usepackage{xcolor}
\usepackage{tikz}
\usepackage{pgfplots}
\usepackage[export]{adjustbox}
\usetikzlibrary{spy}
\usepackage{graphicx}
\usepackage{tabularx}
\usepackage{xcolor}
\usepackage{amssymb,amsmath,amsfonts, amsthm}
\usepackage{tikz}
\usepackage{multirow, hhline}
\usepackage{booktabs}
\usepackage{blindtext}
\usepackage{multicol}
\usepackage{tikzscale}
\usepackage{numprint}
\npdecimalsign{.}
\nprounddigits{3} 
\usepackage{caption}
\usepackage{subcaption}
\usepackage{parskip}
\usepackage{enumitem}
\usepackage{float}
\usepackage{pict2e}
\usepackage{siunitx}
\usepackage[english]{babel}
\usepackage{ulem}
\usepackage{authblk}

\usepackage{empheq}
\usepackage{colortbl}

 \newcommand{\HH}{\mathbf{h}}  
\newcommand{\C}{\mathbb C}

\newcommand{\R}{\mathbb R}

\newcommand{\Rd}{\mathbf{R}}
\newcommand{\CC}{\mathbf{c}}                       
                        
\newcommand{\sd}{\mathbf{s}}                        

\newcommand{\XX}{\mathbf{x}} 
\newcommand{\PP}{\mathbf{p}}

\newcommand{\EE}{\mathbf{e}}                      
\newcommand{\ZZ}{\mathbf{z}}

\newcommand{\YY}{\mathbf{y}}

\newcommand{\Ad}{\mathbf A}                        
\newcommand{\Wd}{\mathbf W} 
\newcommand{\Fd}{\mathbf F} 
                       
\newcommand{\Id}{\mathbf I}                        
\newcommand{\Cd}{\mathbf C}                        
\newcommand{\Ed}{\mathbf E}                        
\newcommand{\Hd}{\mathbf H}

\newcommand{\Au}{\mathbf{A}_I}

\newcommand{\herm}{{\scriptstyle \boldsymbol{\mathsf{H}}}}
\newcommand{\trans}{{\scriptstyle \boldsymbol{\mathsf{T}}}}

\newcommand{\LLambda}{\boldsymbol{\Lambda}} 

\usepackage{mathtools}

\usepackage{mathtools}
\usepackage[percent]{overpic}
\usepackage{algorithm,algcompatible,amsmath}

\DeclareMathOperator*{\argmin}{arg\,min}
\algnewcommand\INPUT{\item[\textbf{Input:}]}%
\algnewcommand\PARAMETER{\item[\textbf{Parameters:}]}%
\algnewcommand\OUTPUT{\item[\textbf{Output:}]}%
\usepackage[
 colorlinks=false, backref=page
]{hyperref}
\renewcommand*{\backref}[1]{}
\renewcommand*{\backrefalt}[4]{%
    \ifcase #1 (Not cited).%
    \or        (Cited on page~#2).%
    \else      (Cited on pages~#2).%
    \fi}

\colorlet{lred}{red!80}
\colorlet{cmix}{blue!80!red} 
\colorlet{lgreen}{green!80}
\colorlet{lblue}{blue!80}

\newcommand{\stkout}[1]{\ifmmode\text{\sout{\ensuremath{#1}}}\else\sout{#1}\fi}
\usepackage{algorithmicx}
\usepackage{comment}

\usepackage{array}

\newcolumntype{L}[1]{>{\raggedright\let\newline\\\arraybackslash\hspace{0pt}}m{#1}}
\newcolumntype{C}[1]{>{\centering\let\newline\\\arraybackslash\hspace{0pt}}m{#1}}
\newcolumntype{R}[1]{>{\raggedleft\let\newline\\\arraybackslash\hspace{0pt}}m{#1}}

\title{Machine Learning for Quantitative MR Image Reconstruction}

\author[1]{Andreas Kofler}
\author[1]{Felix Frederik Zimmermann}
\author[2]{Kostas Papafitsoros}
\affil[1]{Physikalisch-Technische Bundesanstalt (PTB), Braunschweig and Berlin, Germany~~~~~~~~~~}
\affil[2]{School of Mathematical Sciences, Queen Mary University of London, United Kingdom~~~~~\vspace{1em}}
\affil[ ]{\textit {\{andreas.kofler, felix.zimmermann\}@ptb.de}}
\affil[ ]{\textit {k.papafitsoros@qmul.ac.uk}}
\begin{document}

\maketitle

\begin{abstract}
In the last years, the design of image reconstruction methods in the field of quantitative Magnetic Resonance Imaging (qMRI)  has experienced a  paradigm shift. Often, when dealing with (quantitative) MR image reconstruction problems, one is concerned with solving one or a couple of ill-posed inverse problems which require the use of advanced regularization methods. An increasing amount of attention is nowadays put on the development of data-driven methods using Neural Networks (NNs) to learn meaningful prior information without the need to explicitly model hand-crafted priors. In addition,  the available hardware and computational resources nowadays offer the possibility to learn regularization models in a so-called model-aware fashion, which is a unique key feature that distinguishes these models from regularization methods learned in a more classical, model-agnostic manner. Model-aware methods are not only tailored to the considered data, but also to the class of considered imaging problems and nowadays constitute the state-of-the-art in image reconstruction methods. In the following chapter, we provide the reader with an extensive overview of methods that can be employed for (quantitative) MR image reconstruction, also highlighting their advantages and limitations both from a theoretical and computational point of view.
\end{abstract}

\section{Introduction}\label{sec:introduction}
Magnetic Resonance Imaging (MRI) is one the most important medical imaging tools in nowadays clinical practice. MRI allows for the imaging of organs and joints, parallelly exhibiting excellent soft tissue contrast. Unfortunately, the data acquisition process in MRI is inherently slow. In addition, in contrast to other imaging modalities, for example, computed tomography (CT), most MRI scan protocols are not quantitative, i.e.\ the values in the acquired images do not have a physical and/or biophysical correspondence, which represents a challenge for the comparability of images between different scans, scanners, patients or institutions. Quantitative MRI (qMRI) can overcome these limitations by the design of data-acquisition protocols that allow for the quantitative evaluation of biophysical parameters of the imaged matter leading to its subsequent characterization. 

Typically, qMRI is enabled by performing MR measurements for slightly different acquisition configurations, i.e.\ by applying preparation pulses or by varying the MR sequence parameters. By doing so, it is possible to relate the acquired images to several underlying quantitative parameters, such as relaxation times (T1, T2, T2$^\ast$) whose changes were found to be linked to pathological changes in tissue, see e.g.\ \cite{damadian1971tumor, lauterbur1973image} or other references in previous chapters of this book.

This unfortunately results in even longer data acquisition times compared to qualitative MR, which can be addressed by data-undersampling techniques. This, however, yields a series of ill-posed inverse problems which require the use of appropriate regularization methods. In addition, the considered forward model in qMRI is typically considerably more complicated compared to qualitative MRI and the reconstruction requires the use of advanced non-linear reconstruction methods.

In the last years, more and more attention has been put to the development of regularization methods based on data-driven approach, i.e.\ where the regularization approach is learned from the data one wants to reconstruct. This chapter provides an overview of a large variety of methods developed in the last years with a special focus on approaches that make use of Neural Networks.

\subsection{Problem Formulation for Quantitative MR Image Reconstruction}
Let $\XX_{\mathrm{true}} \in \C^N$ denote the vector representation of the (unknown) complex-valued MR image. In the classical (i.e.\ non-quantitative) MRI, often the considered problem is modeled by 
\begin{equation}\label{eq:mri_forward_model}
    \YY = \Ad_J[\Cd] \XX_{\mathrm{true}} + \mathbf{e}, 
\end{equation}
where the linear forward operator is given by 
\begin{eqnarray}\label{eq:measurement_op}
    \Ad_J[\Cd]: \C^N &\longrightarrow & \C^{N_{\mathrm{k}} \cdot N_{\mathrm{c}}}, \\
                \XX &\longmapsto & \big((\Id_{N_c} \otimes \Ed_J) \Cd\big)\, \XX .
\end{eqnarray}
Here the operator $\Cd = [\Cd_1, \ldots, \Cd_{N_c}]^\trans$ with $\Cd_c = \mathrm{diag}(\CC_c), \CC_c \in \C^N, c=1,\ldots, N_c$ denotes the collection of the  $N_c$ coil-sensitivity maps (CSMs), 
which are initially unknown. Furthermore, $\Ed_J: \C^{N}\to \C^{N_{\mathrm{k}}}$ denotes some Fourier encoding operator that samples the coil-weighted images $\Cd_1 \XX_, \ldots, \Cd_{N_c}\XX$ in the Fourier domain (the so-called $k$-space) at arbitrary frequency positions denoted by $J$, i.e.\ either on a Cartesian grid or along other trajectories, e.g.\ radial lines \cite{lauterbur1973image} or spirals \cite{meyer1992fast}. Finally,  $\Id_{N_c}$ denotes the identity operator of dimension $N_c\times N_c$, $\otimes$ denotes the Kronecker product and $\mathbf{e}$ is a random noise component. 

In qMRI, the idea is to collect complementary information about $N_P$ underlying quantitative parameters, which are gathered in a vector $\PP:=[\PP_1, \ldots, \PP_{N_P}]^\trans$ with $\PP \in \mathcal{P}:=\mathcal{P}_1 \times \ldots \times \mathcal{P}_{N_P}$ and $\mathcal{P}_p \in \{\R^N, \C^N\}, p=1,\ldots, N_P$. This information is obtained by repeating the measurement process multiple times with a slightly different configuration, such that the measurements are sensitive with respect to the quantities of interest.

Thus, in qMRI, the considered forward operator can often be described as the composition of a non-linear signal model that encodes the magnetization preparation and a measurement model of the form in \eqref{eq:mri_forward_model} that acquires the data in $k$-space. 
More precisely, let $\mathcal{M}=\{m_1,\ldots, m_{Q}\}$ be a set 
that is used to distinguish between $Q$ different measurements, and $q_{\mathcal{M}}$ a non-linear signal model of the form 
\begin{eqnarray}\label{eq:qmodel}
    q_{\mathcal{M}}: \mathcal{P} &\longrightarrow & \C^{{Q}\,N}, \nonumber\\
                \PP &\longmapsto & [q_{m_1}(\PP), \ldots, q_{m_{Q}}(\PP)]^\trans,
\end{eqnarray}
that maps the $N_P$ different quantitative parameters  to a collection of qualitative images $$[q_{m_1}(\PP), \ldots, q_{m_{Q}}(\PP)]^\trans$$ recorded with different measurement parameters $m_1,\ldots,m_Q.$ For example, the mapping $q_{m_i}$ can be explicitly given as a pre-defined signal model, 
 or can be more complex, e.g.\ the solution of a differential equation\cite{bloch1946nuclear,Wri97,ma2013magnetic,scholand2023quantitative,kim2015review}. However, typically $q_{m_{i}}:\mathcal{P}\to \C^{N}$ acts voxel/pixel-wise, that is, there exists some $\hat{q}_{m_{i}}:\C^{N_{P}}\to \C$ such  that
\begin{equation}\label{q_m_pixel}
q_{m_{i}}(\PP)
=q_{m_{i}}(\PP_{1}, \ldots, \PP_{N_{P}})
=\big ( \hat{q}_{m_{i}}\big(\PP_{1}[j],\ldots,\PP_{N_{P}}[j]\big) \big )_{j=1}^{N}
\end{equation}
where   $\PP_{k}:=\big (\PP_{k}[j] \big )_{j=1}^{N}$ for $k=1, \ldots, N_{P}$.

For example, under some assumptions (see\cite{bloch1946nuclear, Wri97} for more details), two simple examples of signal models that can be used  to describe the evolution of the temporal evolution of the longitudinal and transversal magnetization and yield T1-weighted and T2-weighted images
are given by
\begin{eqnarray}
    \XX_{\mathrm{T1w}} := q_t^{\mathbf{T_1}}(\mathbf{M}_0, \mathbf{T}_1) &= &\mathbf{M}_0\big(\mathbf{1} - \exp\{ -t / \mathbf{T}_1 \}\big),\label{eq:T1_q} \\ 
    \XX_{\mathrm{T2w}} :=  q_t^{\mathbf{T_2}}(\mathbf{M}_0, \mathbf{T}_2) &= & \mathbf{M}_0 \exp\{ -t / \mathbf{T}_2 \},     \label{eq:T2_q}
\end{eqnarray}

where for \eqref{eq:T1_q}, $\PP=[\mathbf{M}_0, \mathbf{T}_1]^\trans$, while for \eqref{eq:T2_q}, $\PP=[\mathbf{M}_0, \mathbf{T}_2]^\trans$ and $\mathbf{M}_0$, $\mathbf{T}_1$, $\mathbf{T}_2$ denote the equilibrium transversal magnetization and the longitudinal and
transverse relaxation relaxation times, respectively.
The equations \eqref{eq:T1_q} and \eqref{eq:T2_q} are solutions of the {\it Bloch equations} \cite{bloch1946nuclear, Wri97}
\begin{equation}\label{eq:bloch}
\begin{aligned}
\frac{\partial \mathbf{M}(t)}{\partial t}
&= \mathbf{M}(t) \times \gamma \mathbf{B}(t) - 
\left( \frac{\mathbf{M}_{z}(t)}{\mathbf{T}_{2}} ,
 \frac{\mathbf{M}_{y}(t)}{\mathbf{T}_{2}} , 
\frac{\mathbf{M}_{z}(t)-\mathbf{M}_{0}}{\mathbf{T}_{1}} 
\right)^{\trans},\\
\mathbf{M}(0)&=(0,0,\mathbf{M}_{0}),
\end{aligned}
\end{equation}

for the particular static magnetic field $\mathbf{B}$ which corresponds to the RF-pulse with $\gamma$ denoting the gyromagnetic ratio.

\begin{figure}[t]
    \centering
    \includegraphics[width=\linewidth]{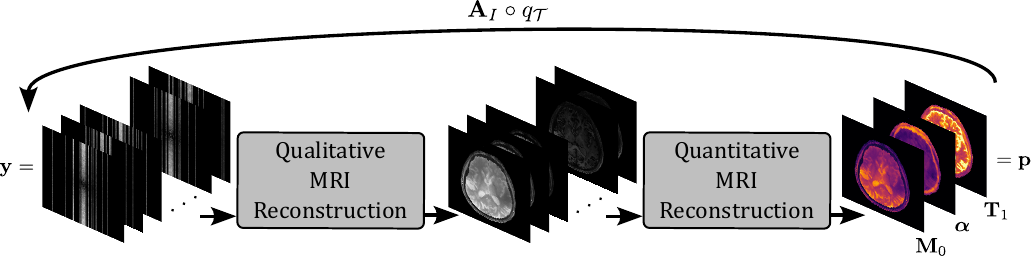}
    \caption{An example for the reconstruction of three different quantitative parameters $\PP=[\mathbf{M}_0, \boldsymbol{\alpha}, \mathbf{T}_1]^\trans$ from undersampled $k$-space data (shown for a single-coil acquisition for simplicity).  Thereby, the set $\mathcal{T}=\{t_1, \ldots, t_{Q}\}$ denotes different inversion times and the quantitative parameter vector $\PP$ contains the equilibrium magnetization $\mathbf{M}_0$, the flip-angle $\boldsymbol{\alpha}$ and the longitudinal relaxation parameter $\mathbf{T}_1$.}
    \label{fig:qmodel}
\end{figure}

Then, the entire forward problem considered in qMRI is of the form 
\begin{equation}\label{eq:qmri_forward_model}
    \YY = \big(\Ad_I[\Cd] \circ q_\mathcal{M} \big)(\PP_{\mathrm{true}}) + \mathbf{e},
\end{equation}
where, by slight abuse of notation,  $\Ad_I[\Cd]$ here denotes
\begin{equation}\label{eq:qmri_forward_model_op}
    \Ad_I[\Cd]:=\mathrm{diag}(\Ad_{I_{m_1}}[\Cd],\, \ldots,\, \Ad_{I_{m_Q}}[\Cd]),
\end{equation}
where each $\Ad_{I_{m_i}}$ is of the form in \eqref{eq:measurement_op} and  where the set of acquired $k$-space points denoted by $I_{m_i}$  potentially varies among different measurements $m_i\in \mathcal{M}$ and
\begin{equation}
    I:= \bigcup_{i=1}^Q I_{m_i}.
\end{equation}

Typically, problems of the form \eqref{eq:mri_forward_model} or \eqref{eq:qmri_forward_model} are ill-posed and require the use of regularization methods to be able to obtain qualitative images or quantitative parameters that can be used for diagnostic purposes. As a first consideration, in the following, we summarize classical variational regularization methods that tackle this ill-posedness for qualitative and quantitative MRI.

\section{Model-based Regularization in  MR Image Reconstruction}\label{sec:model_based}
For the sake of notational simplicity, 
 we neglect the dependence of the operator $\Au[\Cd]$ on the CSMs $\Cd$ in the following, as we are mainly interested in the reconstruction of the qualitative images or quantitative parameters. However, note that obtaining an estimate of $\Cd$ from the $k$-space data already corresponds to solving an inverse problem itself and is often done as a first step before proceeding with the reconstruction of the images/the quantitative parameters, see e.g.\ \cite{uecker2014espirit}. Alternatively, one can jointly reconstruct the images as well as the CSMs using non-linear reconstruction \cite{uecker2008image}.

In the following, we assume to already have an estimate of the CSMs and focus on 
the following coupled inverse problems
\begin{align}
\YY &= \Ad_I \XX + \EE,\label{coupled_ip_1}\\
\XX &= q_{\mathcal{M}}(\PP)+\boldsymbol{\eta}.\label{coupled_ip_2}
\end{align}
Note that while the noise $\EE$ in \eqref{coupled_ip_1} can assumed to be Gaussian, this is not the case for the statistics of the noise component $\boldsymbol{\eta}$ in \eqref{coupled_ip_2}. The reason is that the observed data, i.e.\ the series of qualitative images $\XX$ in \eqref{coupled_ip_2}, is often obtained as a solution to the problem \eqref{coupled_ip_1} using some regularization method. Hence, the qualitative images $\XX$ might exhibit structures and features that are characteristic of the method they were reconstructed with. Therefore, describing the discrepancy between the underlying true and unknown qualitative images and the ones obtained by a prior reconstruction method can be arbitrarily complex and even impossible in practice.

Variational regularization methods have been traditionally employed to deal with the ill-posedness of inverse problems. In the case of \eqref{coupled_ip_1} and \eqref{coupled_ip_2},  the quantities of interest, $\XX$ or $\PP$, are obtained as minimizers of some tailored energy functionals:
\begin{align}
    &\underset{\XX}{\min} \; D(\Ad_I \XX, \YY) + \mathcal{R}(\XX),\label{eq:mri_variational_problem} \\
    &\underset{\PP}{\min} \; D\big((\Ad_I\circ q_{\mathcal{M}})(\PP), \YY) + \mathcal{R}(\PP).\label{eq:qmri_variational_problem}
\end{align}
Here $D(\,\cdot\,,\,\cdot\,)$ denotes a data-discrepancy measure and $\mathcal{R}(\,\cdot\,)$ some regularization term  that imposes desirable properties on $\XX$ or $\PP$ and makes the problem well-posed. In the following, we will always work with the square of the Euclidean norm $D(\,\cdot\,,\,\cdot\,) := \frac{1}{2}\| \cdot \, - \, \cdot\|_2^2$, since it corresponds to the natural metric of choice assuming Gaussian random noise, as a maximum a posteriori estimate (MAP) of the Gaussian likelihood function. 
On the other hand, we will make use a weighted version of it, i.e.\ $\frac{1}{2}\| \cdot \, - \, \cdot\|_{\Wd}^2:= \frac{1}{2}\| \Wd^{1/2} ( \cdot \, - \, \cdot)\|_2^2$ with a positive definite matrix $\Wd$, when sampling the $k$-space data on non-Cartesian grids and employing $k$-space density compensation \cite{pruessmann2001advances}.

In general, directly solving problems 
\eqref{eq:mri_variational_problem} and \eqref{eq:qmri_variational_problem}   can be difficult for different reasons, both from a theoretical as well as computational point of view. For instance, the non-linearity of the mapping $q_{\mathcal{M}}$ can render the minimization \eqref{eq:qmri_variational_problem} non-convex allowing only for computation of local minimizers instead of global ones. Additionally, typical regularization functionals are non-smooth, see Section \ref{sec:reg_qualitative}, raising overall conceptual challenges in the context of algorithmic design. 

To address these issues, there exist several distinct strategies to design regularization and reconstruction methods for the problems in   \eqref{eq:mri_variational_problem} and \eqref{eq:qmri_variational_problem}. We start first with the description of the ones that tackle the simpler qualitative MR image reconstruction problem \eqref{eq:mri_variational_problem} since solution strategies thereof can be employed for solving qMRI problems as in \eqref{eq:qmri_variational_problem} as well.

We note however, that despite the significant progress on non-smooth optimization techniques, practical challenges can emerge even in the simpler problem \eqref{eq:mri_variational_problem}, which is convex for suitable choices of $\mathcal{R}(\XX)$. For example, $\Au$ can be computationally expensive to evaluate. In that case,  its repeated application and/or its adjoint, e.g.\ in the context of primal-dual minimization techniques \cite{chambollepock2016}, should be avoided. Additionally, $\Au$ can be severely ill-conditioned, especially when sampling the $k$-space data along non-Cartesian trajectories.

\subsection{Model-based Regularization for Qualitative MR Image Reconstruction}\label{sec:reg_qualitative}

In qualitative MRI, in view of \eqref{eq:mri_variational_problem} and the choice of the Euclidean norm for data-discrepancy, one is  interested in approximating the solution of \eqref{coupled_ip_1}  by solving
\begin{equation}\label{eq:mri_variational_problem_l2}
    \underset{\XX}{\min}\;  \frac{1}{2}\| \Ad_I \XX - \YY \|_2^2 + \mathcal{R}(\XX).
\end{equation}
Problem \eqref{eq:mri_variational_problem_l2} falls under the remit of ill-posed linear inverse problems, an area with a vast literature \cite{benning_burger_2018, engl1996regularization}. Here the regularization term $\mathcal{R}(\,\cdot\,)$ encodes some a priori information that is imposed on the reconstructed image. This information is typically related to some function regularity (e.g.\ smoothness or discontinuity), sparsity in some basis, or in general related to some specific structure. One of the classical choices is that of the total variation (TV) 
\begin{equation}\label{def:TV}
\mathcal{R}(\XX)=\lambda \|\nabla \XX\|_{1}=:\lambda\mathrm{TV}(\XX),
\end{equation}
originating from the early 90s for image denoising, and later applied to MRI reconstruction \cite{rudin1992nonlinear, Lustig_2007, ChambolleLions, block2007undersampled}. Here $\|\cdot\|_{1}$ denotes the discrete $\ell_{1}$-norm, with $\|r\|_{1}:=\sum_{i}|r_{i}|$, $r\in\mathbb{R}^{s N}$, $r_{i}\in\R^{s}$, where $s$ denotes the number of directions for which the partial derivatives are computed, i.e.\ $s=2$ or $3$ for 2D and 3D images respectively.
The scalar parameter $\lambda>0$ balances the effect of the fidelity term and TV in \eqref{eq:mri_variational_problem_l2}. The literature on TV is itself quite extensive \cite{Caselles2015, scherzer2009variational, Burger2013}. We mention here that the use of $\ell_{1}$-norm enforces sparsity in the gradient domain, promoting piecewise constant reconstructions \cite{Jalalzai2015jmiv, poon_TV_geometric, ring2000structural}. This has the advantage of edge preservation but it also leads to undesirable artifacts (staircasing effect), which in the case of medical imaging can interfere with the diagnostic procedure.
Towards suppressing the later effect, regularization functionals incorporating higher order derivatives have been suggested with a prominent example that of Total Generalized Variation (TGV) \cite{TGV}. TGV regularization promotes piecewise affine structures in the reconstruction
and has been extensively considered for both static and dynamic MRI reconstruction \cite{knoll2011second, TGV_dynamic_MRI}.

 More sophisticated regularizers can also be considered. For example, infimal convolution-type regularization promotes images that can be decomposed into components of different structures (e.g.\ smooth and piecewise constant) \cite{journal_tvlp, ChambolleLions} and highly directional regularizers can enhance anisotropic features, see \cite{Parisotto_2020} and the references therein.  Regularizers that enforce sparsity in some basis, e.g. wavelet, shearlets \cite{Donoho_1994, Chang_2000, Kutyniok2018} have also been popular in MRI reconstruction.

In all of these approaches, the regularization parameter(s)  have to be chosen carefully. For instance, in the case of TV, larger $\lambda$ lead to over-regularization (smoothing) and a loss of details while smaller $\lambda$ result in poor regularization (fitting to noisy data) and insufficient suppression of artifacts. Crucially,  the regularization strength encoded by these parameters should not be applied uniformly across the image domain, since images consist of both homogeneous areas and fine-scale details. This gives rise to regularization functionals with spatially dependent regularization maps, which in the case of TV would read
\begin{equation}\label{def:TV_spatial}
\mathcal{R}(\XX)= \|\LLambda\nabla \XX\|_{1}.
\end{equation}
Here, $\LLambda\in \R_{+}^{sN}$ is a diagonal matrix denoting a spatially varying  (pixel/voxel dependent) regularization parameter. This type of regularization has achieved significant attention in the literature with regard to theoretical investigations \cite{jalalzai2014discontinuities, papafitsoros_spatial}, and extensions to other regularizers beyond TV \cite{bilevelconvex, bilevelTGV}. An important challenge is the automatic selection of the map $\LLambda$ which due to its large number of components becomes a non-trivial task. Towards that, bilevel optimization techniques have been employed during the last years, aiming to compute a data-adapted $\LLambda$ via the minimization of an upper-level objective $l$, see \cite{DelosReyes2021, Fessler, bilevel_handbook, bilevellearning} and also \cite{Pragliola_SIAM_review} for some additional approaches. These have the following general formulation: 
\begin{equation}\label{general_bilevel}
\left \{
\begin{aligned}
&\min_{\boldsymbol{\Lambda}}\; \frac{1}{|\mathcal{D}|}\sum_{i=1}^{N_{\mathrm{train}}}l(\XX^{i}(\boldsymbol{\Lambda}),\XX_{\mathrm{true}}^{i})\\
&\text{subject to }\;\;\XX^{i}(\boldsymbol{\Lambda})=\argmin_{\XX'}\;\frac{1}{2} \| \Ad_I \XX' - \YY^{i}\|_2^2 +  \| \boldsymbol{\Lambda}\nabla \XX' \|_1, \quad i=1,\ldots,N_{\mathrm{train}}.
\end{aligned}\right.
\end{equation}
Here, $\mathcal{D}:=\{(\YY^{i}, \XX_{\mathrm{true}}^{i})_{i=1}^{N_{\mathrm{train}}}\}$ are $N_{\mathrm{train}}$ pairs of measured data  and corresponding ground truth images. For instance, in the case where $l(\XX_{1},\XX_{2})=l_{\mathrm{PSNR}}(\XX_{1},\XX_{2}):=\|\XX_{1}-\XX_{2}\|_{2}^{2}$, the bilevel problem \eqref{general_bilevel} aims to compute the parameters $\boldsymbol{\Lambda}$ which are ``on the average the best ones" (i.e.\ PSNR-maximizing), for the given $N_{\mathrm{train}}$ data pairs 
. We note however that there are 
approaches for which the upper-level objective does not involve the ground truth, see \cite{bilevel_handbook}. Another recent approach,  also mentioned later in Section \ref{sec:data_driven} employs neural networks to compute such parameters also for dynamic image reconstruction problems including dynamic cardiac MRI and qMRI \cite{Kofler_SIAM_Imaging_2023}. In this case, the regularization parameter-maps are spatio-temporally varying.

The use of the above hand-crafted regularization approaches for qualitative MRI has several advantages:
\begin{enumerate}
    \item {\it Intepretability}: The reconstructed image is fully interpretable as a solution of the minimization problem \eqref{eq:mri_variational_problem_l2}. It also possesses an a priori known and expected structure which is imposed by the regularizer.
    \item {\it Guarantees}: There exist nowadays numerous algorithms for solving \eqref{eq:mri_variational_problem_l2} with  convergence guarantees  \cite{chambollepock2016}. Additionally, in the limit of vanishing noise $\EE\to 0$ and the regularization parameters converging to zero at an appropriate rate, the solutions of \eqref{eq:mri_variational_problem_l2} converge to an image that is data-consistent, i.e.\, $\YY = \Ad_I \XX$. 
    \item {\it Control of artifacts}: Potential artifacts, e.g. staircasing effect, are well-understood and can be taken into account during the diagnostic procedure. They can be further controlled with regularization of spatially varying strength.
    \item {\it Amenable to rigorous mathematical analysis}: Problem \eqref{eq:mri_variational_problem_l2}, TV/TGV/wavelet based regularization functionals as well as bilevel optimization problems like \eqref{general_bilevel} can all be studied in the infinite-dimensional setting. There, the notion of discontinuity/sharp edges can be rigorously defined, and theoretical analysis provides information about whether the artifacts are due to the model or due to the numerical discretization, (e.g.\ numerical diffusion \cite{functional_issues}). 
\end{enumerate}

On the other hand, these approaches have also two main disadvantages:

\begin{enumerate}
    \item {\it Non-realistic image priors}: The imposed structure and prior information of these regularizers (e.g. piecewise constant/affine) might not reflect features and structures of real-world MR images well. As a result, these methods do not typically produce state-of-the-art reconstructions, especially when compared to data-driven methods mentioned in Sections \ref{sec:data_driven} and \ref{sec:nns_for_mri}.
    \item {\it Computational cost}: Solving bilevel optimization problems requires repeatedly solving a minimization problem, which can be computationally demanding and inefficient. As most of these reconstruction methods involve the use of iterative schemes to obtain the solution, depending on the considered forward model, this could amount to several hours or even days for large-scale  problems, e.g. dynamic problems and/or 3D imaging problems. 
    \end{enumerate}

\subsection{Model-based Regularization for Quantitative MR Image Reconstruction}

In the following, we present two main approaches for estimating $\PP$ depending on whether the problems \eqref{coupled_ip_1} and \eqref{coupled_ip_2} are solved sequentially (decoupling approaches) or they are coupled. In the following, we denote by  $q_{\mathcal{M}}$  the non-linear signal model introduced in \eqref{eq:qmodel}  used to collect complementary information about the pixel-wise signal evolution for different measurements set-ups defined by the set $\mathcal{M}:=\{m_1, \ldots, m_Q\}$.

\subsubsection{Decoupling Approaches}

An often pursued strategy is to decouple  $q_{\mathcal{M}}$ from the 
linear data-acquisition operator $\Au$ and to recast the problem in the following form 
\begin{empheq}[left=\empheqlbrace]{alignat=1} 
&\underset{\PP}{\min}\;  \frac{1}{2}\| q_{\mathcal{M}}(\PP) - \XX \|_2^2, \label{eq:qmri_variational_problem_splitted1}\\ 
&\text{where}\;\; \XX = \underset{\XX^\prime}{\argmin} \; \frac{1}{2}\| \Ad_I \XX^\prime - \YY \|_2^2 + \mathcal{R}(\XX^\prime). \label{eq:qmri_variational_problem_splitted2} 
\end{empheq}
By doing so, all regularization effort is put in the reconstruction of the series of qualitative images $\XX:=[\XX_1, \ldots, \XX_Q]^\trans$, resulting in $Q$ distinct qualitative MRI reconstruction problems, and in principle, every method discussed in Subsection \ref{sec:reg_qualitative} can be employed for that purpose. Notice that the problems \eqref{eq:qmri_variational_problem_splitted1}--\eqref{eq:qmri_variational_problem_splitted2} can be solved sequentially: After solving \eqref{eq:qmri_variational_problem_splitted2} and having obtained the qualitative images, the set of quantitative parameters $\PP$ can be estimated using a non-linear optimization routine, e.g.\ the limited memory Broyden–Fletcher–Goldfarb–Shanno (L-BFGS) algorithm \cite{liu1989limited, nocedal1999numerical} or the Levenberg–Marquardt algorithm \cite{marquardt1963algorithm}. 

The {\it Magnetic Resonance Fingerprinting} (MRF) method \cite{ma2013magnetic}, considered in its original form, follows closely the form of \eqref{eq:qmri_variational_problem_splitted1}--\eqref{eq:qmri_variational_problem_splitted2}.
At a first step, a series of qualitative MR images resulting from a long and complex sequence of RF pulses is reconstructed via a least square solution, i.e.\ $\mathcal{R}(\XX)\equiv 0$ in \eqref{eq:qmri_variational_problem_splitted2}.
The problem \eqref{eq:qmri_variational_problem_splitted1} is then solved by pre-computing, in an offline phase, responses of the model $q_{\mathcal{M}}$ for $N_{\mathrm{dic}}$
predefined values for $\PP$ on the pixel/voxel level $\hat{\PP}^{\ell}=[\hat{\PP}_{1}^{\ell}, \ldots, \hat{\PP}_{N_{P}}^{\ell}]\in \C^{N_{P}}$, $\ell=1,\ldots, N_{\mathrm{dic}}$ forming a dictionary (the {\it fingerprints}) $(\hat{q}_{\mathcal{M}}(\hat{\PP}^{\ell}))_{\ell=1}^{N_{\mathrm{dic}}}$, with every $\hat{q}_{\mathcal{M}}(\hat{\PP}^{\ell})\in \C^{N_{Q}}$. Here $\hat{q}_{\mathcal{M}}(\hat{\PP}^{\ell}):=(\hat{q}_{m_{1}}(\hat{\PP}^{\ell}), \ldots, \hat{q}_{m_{N_{Q}}}(\hat{\PP}^{\ell}))$ recall \eqref{q_m_pixel}.
This results to 
\begin{empheq}[left=\empheqlbrace]{alignat=1} 
&\PP[j]=\argmin_{\hat{\PP}\in \big\{\hat{\PP}^{\ell}\big\}_{\ell=1}^{N_{\mathrm{dic}}}}\; \frac{1}{2} \big\|\hat{q}_{\mathcal{M}}(\hat{\PP})-\XX[j]\big\|_{2}^{2}, \quad j=1,\ldots, N\label{eq:qmri_variational_problem_splitted1_MRF}\\
&\text{where}\;\; \XX = \underset{\XX^\prime}{\argmin} \; \frac{1}{2}\| \Ad_I \XX^\prime - \YY \|_2^2 . \label{eq:qmri_variational_problem_splitted2_MRF} 
\end{empheq}
Here, for all $j=1,\ldots,N$ we have $\XX[j]\in \C^{N_{Q}}$, recall \eqref{eq:qmri_forward_model} and \eqref{eq:qmri_forward_model_op}. One of the main disadvantages of the original MRF method is the large storage requirements for the dictionary which needs to be rich enough (i.e.\ $N_{\mathrm{dic}}$ should be large), with the dictionary matching itself being computationally demanding. Additionally, the poor quality of the qualitative reconstructions \eqref{eq:qmri_variational_problem_splitted2_MRF} can have an adverse effect on the inferred $\PP$. Several extensions of MRF have been considered aiming to improve computational times and/or reconstruction quality for $\PP$. For instance, iterative solution schemes for \eqref{eq:qmri_variational_problem_splitted2_MRF} with $\mathcal{R}\equiv 0$ as well as with wavelet spatial regularization were considered in \cite{Davies_2014},   also involving a projection (matching) of each qualitative image onto the dictionary at every step. The use of fast searching techniques for this projection \cite{Golbabaee_2020} has also been considered, involving low-rank compression of the dictionary \cite{mcgivney2014svd} or the image itself \cite{Mazor_2018}. We refer to the next sections regarding extensions of MRF that employ neural networks.

\begin{figure}[t!]
    \centering
\includegraphics[width=\linewidth]{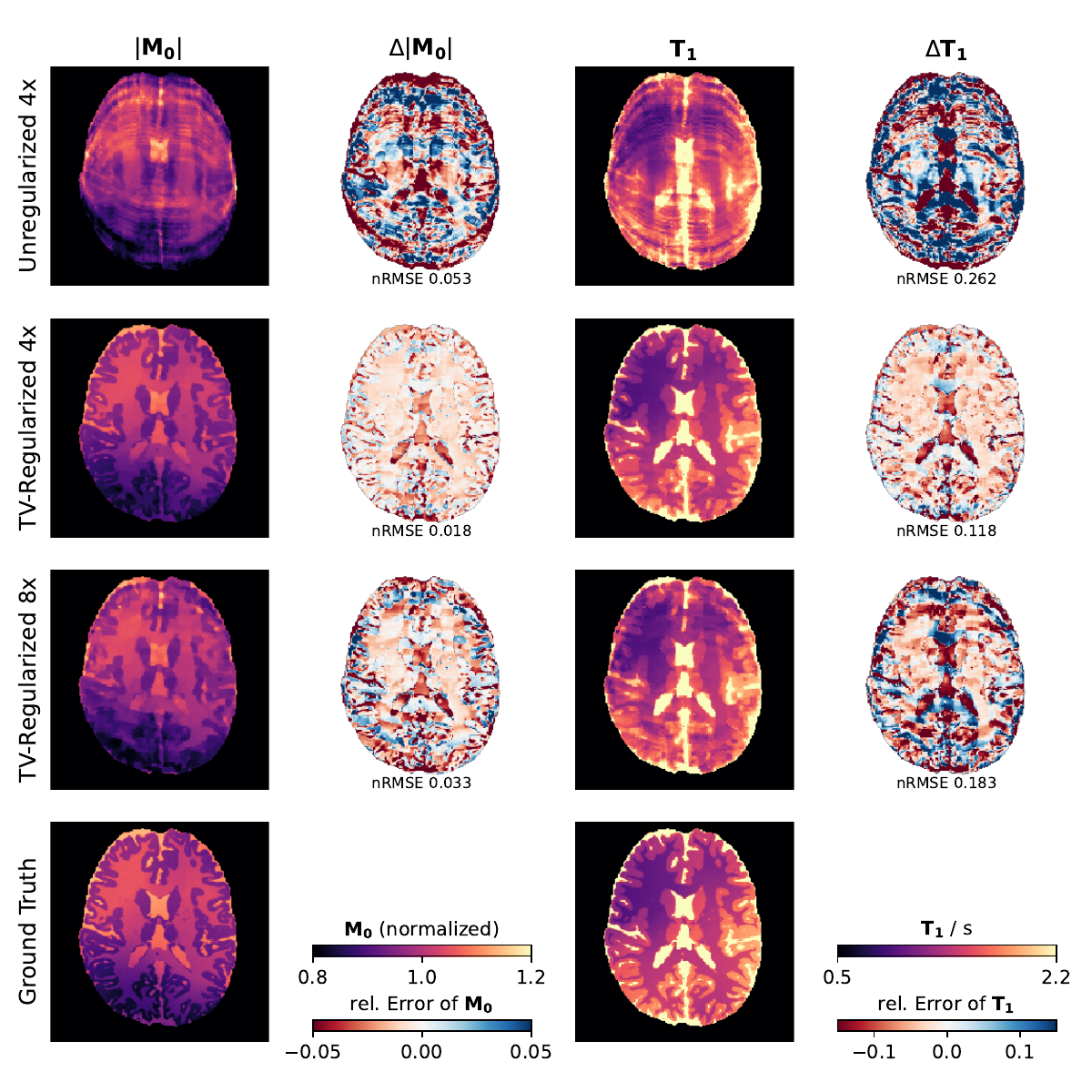}
    \caption{A comparison of examples for the reconstruction of $\mathbf{M}_0$ and $\mathbf{T}_1$ from simulated 4-fold and 8-fold Cartesian undersampled saturation recovery data. Without any regularization, the results obtained at 4-fold undersampling (top row) are severely degraded by artifacts. TV-regularization, i.e.\. setting $\mathcal{R}(\XX)$ according to \eqref{def:TV} in \eqref{eq:qmri_variational_problem_splitted2}, improves the results visibly at 4-fold (second row), but fails at 8-fold undersampling (third row). This is the setting that is used to demonstrate different NN-based approaches in Figure \ref{fig:model_aware_vs_model_unaware_vs_pinqi_qmri}. The last row shows the ground-truth target labels of $\mathbf{M}_0$ and $\mathbf{T}_1$ used in the simulation (generated from the BrainWeb dataset\cite{aubert2006twenty}).
      }
    \label{fig:tv+results}
\end{figure}

Conceptually, decoupling strategies can be attractive from a computational point of view, because the original problem \eqref{eq:qmri_variational_problem_l2} is divided into two problems: the first being a relatively large problem but where the forward operator is 
linear and well-studied in the literature and the second being a pixel-wise non-linear operator. This strategy avoids the need to repeatedly apply full operator $\Au\circ q_{\mathcal{M}}$, as it would, for example, be the case for gradient descent-type methods that involve line-searches.
An example of results obtained by a decoupled approach with TV-regularization for the qualitative images is shown in Figure \ref{fig:tv+results}. The effect of the regularization is visible by the improved reconstruction at 4-fold acceleration compared to an unregularized approach.

However, at the same, the dimensionality of the qualitative image reconstruction problem is rather large, as it involves the reconstruction of $Q$ qualitative images. This aspect can be especially challenging depending on the problem formulation of the qualitative image reconstruction problem. In addition, this strategy inherently introduces a bias that is dependent on the regularization method of choice used in \eqref{eq:qmri_variational_problem_splitted2}. This also means that systematic errors introduced in the qualitative images due to the employed regularization cannot be counteracted in the non-linear fitting routine and are potential sources for inaccuracies in the quantification of the parameters contained in $\PP$.

There exists a large variety of further approaches that follow the discussed strategy to tackle image reconstruction and quantitative parameter mapping as two separate problems. For example, in \cite{petzschner2011fast}, problem \eqref{eq:qmri_variational_problem_splitted2} is solved using $k$-$t$-PCA \cite{pedersen2009k} for a T1- and T2-mapping, and in \cite{becker2019simultaneous}  $k$-$t$-SENSE \cite{tsao2003k} was used to reconstruct cardiac MR images for T1-mapping. As mentioned already, in \cite{Kofler_SIAM_Imaging_2023} problem \eqref{eq:qmri_variational_problem_splitted2} was solved with a spatially varying TV regularization where the regularization parameter-maps are inferred by a neural network.

\subsubsection{Coupling Approaches}

Instead of  reconstructing a series of regularized qualitative images first and subsequently performing a non-linear fit to obtain $\PP$, it is also possible to directly solve the least squares problem
\begin{equation}\label{eq:qmri_fitting_problem_l2}
    \underset{\PP}{\min}\;  \frac{1}{2}\| (\Ad_I \circ q_{\mathcal{M}})(\PP) - \YY \|_2^2,
\end{equation}
using an iterative scheme that imposes implicit regularization, e.g.\ a projected Levenberg-Marquardt method \cite{Hanke_1997, Dong_2019}:
\begin{align}
\tilde{\YY}_{k}&=\YY-(\Ad_I\circ q_{\mathcal{M}})(\PP_{k}),\label{LM1}\\ 
\HH_{k}&=\argmin_{\HH}\; \frac{1}{2}\,\|(\Ad_I\circ q_{\mathcal{M}})'(\PP_{k})\HH-\tilde{\YY}_{k}\|_{2}^{2} + \frac{\lambda_{k}}{2}\|\HH\|_{2}^{2},\label{LM2}\\
\PP_{k+1}&=P_{C_{\mathrm{ad}}}(\PP_{k}+\HH_{k}), \label{LM3}
\end{align}
where $\{\lambda_{k}\}_{k}$ is a sequence of positive parameters that converges to zero.
This is an iterative Tikhonov regularization method for solving the non-linear equation $(\Ad_I\circ q_{\mathcal{M}})(\PP)=\YY$ and requires some early stopping based on a discrepancy principle \cite{Hanke_1997}. The operation $P_{C_{\mathrm{ad}}}$ in \eqref{LM3} denotes a projection to  an admissible domain $C_{\mathrm{ad}}$ for the values of $\PP$, encapsulating realistic biophysical values.
The sequence  $\{\lambda_{k}\}_{k}$, as well as the initialization $\PP_{0}$ need to be chosen carefully to achieve good reconstructions. In \cite{Dong_2019}, it was shown that a {\it warmly} initialized iterative regularization approach as in \eqref{LM1}--\eqref{LM3} produces improved results compared to the MRF ones \eqref{eq:qmri_variational_problem_splitted1_MRF}--\eqref{eq:qmri_variational_problem_splitted2_MRF}.

Another approach is to formulate the problem 
\begin{equation}\label{eq:qmri_variational_problem_l2}
    \underset{\PP}{\min}\;  \frac{1}{2}\| (\Ad_I \circ q_{\mathcal{M}})(\PP) - \YY \|_2^2  + \mathcal{R}(\PP)
\end{equation}
and directly impose certain properties on the sought quantitative parameters via the regularization term $\mathcal{R}$. Essentially, this corresponds to \eqref{eq:qmri_variational_problem_splitted1}--\eqref{eq:qmri_variational_problem_splitted2} where the minimization \eqref{eq:qmri_variational_problem_splitted1} is replaced by the hard constraint $q_{\mathcal{M}} (\PP) = \XX$ which is then inserted in  \eqref{eq:qmri_variational_problem_splitted2}.
Prominent examples for $\mathcal{R}(\PP)$ include the ones previously discussed in Subsection \ref{sec:reg_qualitative}, e.g.\  $\mathcal{R}(\PP)=\mathrm{TGV}(\PP)$ \cite{wang2018model}, $\mathcal{R}(\PP) = \| \nabla \PP\|_2^2$ \cite{olafsson2008fast, Dong_2022}, as well as $\mathcal{R}(\PP)=\| \nabla \Fd \PP\|_2^2$ with $\Fd$ being the 2D FFT operator \cite{block2009model} or $\mathcal{R}(\PP) = \sum_{p=1}^{P} \| \Wd \PP_p \|_1 $ with $\Wd$ being a Wavelet transform \cite{wang2018model}.
As the entire considered forward model is non-linear, these type of approaches require the use of advanced strategies for iteratively minimizing the objective function \eqref{eq:qmri_variational_problem_l2}, e.g.\ see \cite{nocedal1999numerical}.

\subsection{ Linear Subspace Methods}

Another type of approach widely used to overcome the non-linearity of the signal model $q_{\mathcal{M}}$ are linear sub-space methods, see e.g.\ \cite{huang2012t2, pfister2019simultaneous, tamir2017t2}. In these, one constructs a temporal basis of the signal evolution of the qualitative images generated according to $q_{\mathcal{M}}$. Thereby, based on the observation that the qualitative images generated by $q_{\mathcal{M}}$ have low-dimensional structure and are thus well-representable by a linear combination of only a few basis functions, one substitutes
\begin{equation}\label{eq:subspace_approx}
    q_{\mathcal{M}}(\PP) \approx \sum_{j=^1}^K  \mathbf{\Psi}_k \mathbf{s}_k,
\end{equation}
where $\{\mathbf{\Psi}\}_{k=1}^K$ are the first $K$ basis functions and $\{\mathbf{s}\}_{k=1}^K$ the corresponding coeffcient-images. By inserting the representation \eqref{eq:subspace_approx} in \eqref{eq:qmri_variational_problem_l2}, the problem becomes linear and one can employ any of the previously discussed regularization methods.

\section{Data-Driven Dictionary-Based  Regularization Methods}\label{sec:data_driven}

The urge to design data-driven regularization methods that are adapted to the data under consideration stems from the observed limitations of hand-crafted image priors such as the previously discussed examples, e.g.\ TV, TGV, Wavelets, Shearlets, etc. In general, the idea is to replace the term $\mathcal{R}$ by a term that depends on parameters $\Theta$ to be learned from data, i.e.\ $\mathcal{R}_{\Theta}$.
Here we outline some methods that are data-driven but still hand-crafted, focusing mainly on dictionary learning, leaving the focus on the more modern neural network-based methods for the later sections. Again, we start with the presentation of the methods  
for qualitative MRI and continue with their adaptation for quantitative MRI. 

\subsection{Data-Driven Dictionary-Based Regularization Methods for Qualitative MRI}\label{sec:data_driven_mri}

A prominent example that uses sparsity as regularization is dictionary learning \cite{rubinstein2010dictionaries}. A dictionary is referred to as a (typically overdetermined) basis of $K$ elements used to approximate signals with using only $S$ of its $K$ elements, where $S\ll K$, see \cite{rish2014sparse} for an extensive introduction to sparse modeling. The dictionary learning problem, i.e. to learn a dictionary that is suitable for the sparse representation of $N_{\mathrm{train}}$ $d$-dimensional signals $\{\ZZ_j\}_{i=1}^{N_{\mathrm{train}}}$, can be formulated as 
\begin{equation}\label{eq:dico_learning_problem_L0}
\underset{\mathbf{\Psi} \in \mathcal{D}_{d,K}, \, \{\boldsymbol{\gamma}_j\}_j}{\mathrm{min}} \frac{1}{2}\sum_{j=1}^{N_{\mathrm{train}}} \| \ZZ_j - \mathbf{\Psi} \boldsymbol{\gamma}_j\|_2^2  \quad \text{ such that } \quad   \forall j: \| \boldsymbol{\gamma}_j\|_0 \leq S,  
\end{equation}
where $\mathcal{D}_{d,K}:=\{ \mathbf{\Psi} \in \mathbb{R}^{d\times K} :  \|\boldsymbol{\psi}_k \|_2 = 1, k=1,\ldots,K\}$ denotes the set of admissible dictionaries, that is all $d$-dimensional dictionaries with $K$ atoms with unit-norm.

In the context of MR image reconstruction, the signals to be sparsely approximated by the learned dictionary  $\mathbf{\Psi}$ are patches that are extracted from the image. Thus, the underlying regularizing assumption is that patches of MR images have an inherently low-dimensional structure and that artifacts and noise cannot be well-approximated by the learned dictionary.

For qualitative MR image reconstruction, dictionary learning has been extensively applied to different organs, e.g.\ cardiac MRI \cite{caballero2014dictionary, wang2013compressed, pali2021adaptive} or brain MRI \cite{ravishankar2010mr, song2014reconstruction,song2019coupled}. Thereby, there exist two main different strategies to employ dictionary learning for regularization. In the first, the dictionary is trained beforehand on a set of training-patches and used for inference, while in the second, the dictionary is learned during the reconstruction and thus is also adapted to the image that is being reconstructed. The corresponding problem formulations are 
\begin{equation}\label{eq:dlmri_psi_fixed}
\underset{\XX, \{\boldsymbol{\gamma}_j\}_j}{\mathrm{min}} \frac{1}{2}\| \Au  \XX -\YY\|_2^2 + \frac{\lambda}{2} \sum_{j=1}^{N_{\mathrm{patches}}} \| \Rd_j \XX - \mathbf{\Psi} \boldsymbol{\gamma}_j \|_2^2 \quad \text{such that} \quad  \forall j: \| \boldsymbol{\gamma}_j\|_0 \leq S,
\end{equation}
and
\begin{equation}\label{eq:dlmri_psi_ada_learned}
\underset{\XX, \mathbf{\Psi} \in \mathcal{D}_{d,K}, \{\boldsymbol{\gamma}_j\}_j}{\mathrm{min}} \frac{1}{2}\| \Au  \XX -\YY\|_2^2 + \frac{\lambda}{2} \sum_{j=1}^{N_{\mathrm{patches}}} \| \Rd_j \XX - \mathbf{\Psi} \boldsymbol{\gamma}_j \|_2^2 \quad \text{such that} \quad   \forall j: \| \boldsymbol{\gamma}_j\|_0 \leq S,
\end{equation}
where $\Rd_j$ extracts the $j$-th patch from the image. Problems \eqref{eq:dlmri_psi_fixed} and \eqref{eq:dlmri_psi_ada_learned} are most commonly solved using alternating minimization. Thereby, the sub-problem of \eqref{eq:dlmri_psi_ada_learned} with respect to the dictionary $\mathbf{\Psi}$ and the set of sparse codes $\{\boldsymbol{\gamma}_j\}_j$, exactly corresponds to problem \eqref{eq:dico_learning_problem_L0}. Reconstruction methods derived from problem formulations similar to \eqref{eq:dlmri_psi_ada_learned} are sometimes also referred to as \textit{blind} Compressed Sensing approaches because the dictionary is not known beforehand, but instead learned during the reconstruction, see e.g.\ \cite{ravishankar2010mr, wang2013compressed, caballero2014dictionary}. Importantly, we note that in \eqref{eq:dlmri_psi_ada_learned}, the regularization method is an optimization problem itself that depends on the set of trainable parameters, i.e.\ the entries of the dictionary.

\begin{figure}
    \centering
    \includegraphics[width=\linewidth]{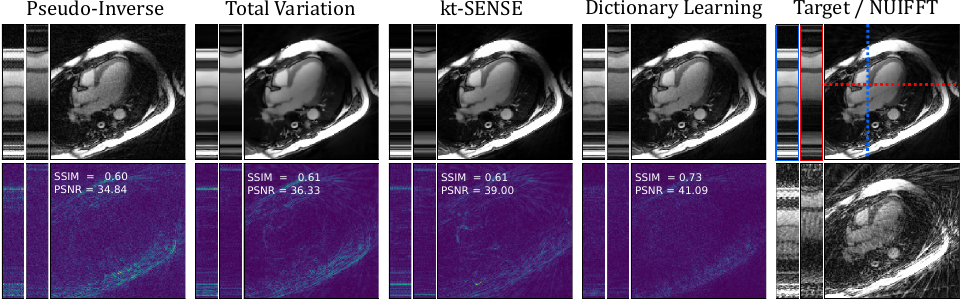}
    \caption{An example of cardiac cine MR images reconstructed from undersampled radially acquired $k$-space data. From left to right: approximation of the Moore-Penrose pseudo-inverse, TV-minimization \cite{block2007undersampled}, $kt$-SENSE \cite{tsao2003k} and a dictionary learning-based approach \cite{pali2021adaptive}. Additionally, the target image from which the $k$-space data was retrospectively simulated as well as the initial reconstruction, given by applying the adjoint operator to the density-compensated $k$-space data, are also shown. The data-driven dictionary learning approach improves the images compared to the other shown methods both in terms of PSNR as well as SSIM.}
    \label{fig:moore_vs_tv_kt_sense_vs_dl}
\end{figure}

Figure \ref{fig:moore_vs_tv_kt_sense_vs_dl} shows an example of dynamic cardiac MR images reconstructed with different approaches, i.e.\ with the approximation of the Moore-Penrose pseudo-inverse obtained by approximately solving the normal equations in \eqref{eq:normal_eqs}, a TV-minimization approach \cite{block2007undersampled}, $kt$-SENSE \cite{tsao2003k} and a dictionary learning-based approach with adaptive dictionary-size and sparsity level choice \cite{pali2021adaptive} at an acceleration factor of approximately nine. The improvements of the data-driven dictionary learning method over the hand-crafted image priors used in the TV and $kt$-SENSE reconstruction is evident from the point-wise error images. In addition, the structural similarity index measures (SSIM) \cite{wang2004image} improves substantially.

Despite the promising results that these types of approaches deliver, there remain some major difficulties. Learning the dictionary can be time-consuming, which is an important factor to take into account when the dictionary is adaptively learned during the reconstruction as in \eqref{eq:dlmri_psi_ada_learned}. The most prominent dictionary learning algorithm is $K$-SVD \cite{aharon2006k}, which updates the dictionary atoms using $K$ singular value decompositions (SVDs) and updates the sparse codes by orthogonal matching pursuit \cite{pati1993orthogonal}. The time required for the learning dictionary with $K$-SVD additionally depends on the choice of $S$ and $K$ and typically, larger $S$ and $K$ yield longer computational times.
To that end, alternative dictionary algorithms, e.g.\ the iterative thresholding and $K$ residual means (ITKrM) \cite{schnass2018convergence}, can be employed to reduce the time required for learning the dictionary.

However, the real computational bottleneck of these methods lies in the need to repeatedly perform the sparse approximation of all $N_{\mathrm{patches}}$ image patches when solving \eqref{eq:dlmri_psi_fixed} or \eqref{eq:dlmri_psi_ada_learned} for $\{\boldsymbol{\gamma}_j\}_j$. Employing sparsity- and dictionary-size-adaptive sparse approximation and dictionary learning algorithms \cite{pali2023dictionary} was shown to be a promising alternative to further reduce reconstruction times, see \cite{pali2021adaptive} for a comparison of the overall computational time required by the algorithms used to solve the sub-problems of \eqref{eq:dlmri_psi_ada_learned}. 
Last, the employed  $\ell_0$-pseudo-norm $\|\, \cdot, \|_0$ is not convex. Convex-relaxation methods that substitute the $\ell_0$-pseudo-norm by the convex $\ell_1$-norm address this issue by replacing the non-convex problem with a convex one
have been considered as well \cite{kofler2022nndlmri}.

To further overcome the computational difficulties related to these patch-based approaches,  the convolutional pendant for sparse approximation, i.e.\ convolutional dictionary learning can be employed. Convolutional dictionary learning, instead of being a local model as the patch-based dictionary learning, is a global model and assumes that an entire image can be well-approximated by the convolution of dedicated filters and sparse images \cite{wohlberg2015efficient}. Given a dataset of $N_{\mathrm{train}}$ target images $\{\ZZ_j\}_{j=1}^{N_{\mathrm{train}}}$, the corresponding convolutional dictionary learning problem using the $\ell_1$-norm for measuring sparsity is typically posed as 
\begin{equation}\label{eq:cd_learning_problem_L1}
\underset{ \{ d_k\}_k \in \mathcal{D}^{\ast}_{d,K}, \, \{\sd_{k,j}\}_{k,j}}{\mathrm{min}}\, \frac{1}{2}\sum_{j=1}^{N_{\mathrm{train}}} \Bigg( \| \ZZ_j - \sum_{k=1}^K d_k \ast \sd_{k,j}\|_2^2  + \alpha \sum_{k=1}^K\| \sd_{k,j}\|_1\Bigg),  
\end{equation}
where $\mathcal{D}^{\ast}_{d,K}:= \big\{ \{d_k\}_{k=1}^K : \forall k: \|d_k \|_2 = 1 \big\}$ denotes the set of $K$ filters with unit-norm and $\alpha>0$. When used for MR reconstruction, it can be employed similarly as the previously introduced patch-based dictionary learning, by formulating the problems
\begin{equation}\label{eq:cdlmri_psi_fixed}
\underset{\XX, \{\sd_k\}_k}{\mathrm{min}} \frac{1}{2}\| \Au  \XX -\YY\|_2^2 + \frac{\lambda}{2} \| \XX - \sum_{k=1}^K d_k \ast \sd_{k}\|_2^2  + \alpha \sum_{k=1}^K\| \sd_{k}\|_1,
\end{equation}
and
\begin{equation}\label{eq:cdlmri_psi_ada_learned}
\underset{\XX, \{ d_k\}_k \in \mathcal{D}^{\ast}_{d,K}, \, \{\sd_k\}_k}{\mathrm{min}} \frac{1}{2}\| \Au  \XX -\YY\|_2^2 + \frac{\lambda}{2} \| \XX - \sum_{k=1}^K d_k \ast \sd_{k}\|_2^2  + \alpha \sum_{k=1}^K\| \sd_{k}\|_1'
\end{equation}
Again, in \eqref{eq:cdlmri_psi_fixed} the set of convolutional dictionary filters is fixed and obtained beforehand by solving \eqref{eq:cd_learning_problem_L1} on a set of target images, while in \eqref{eq:cdlmri_psi_ada_learned}, the convolutional dictionary is adaptively learned during the reconstruction. We refer to \cite{garcia2018convolutional} for an extensive review on methods for convolutional dictionary learning and to \cite{quan2016compressed1, nguyen2018compressed, kofler2022cdl} for its application to qualitative MR reconstruction.

Last, we note that dictionary learning, regardless of patch-based or convolutional, only corresponds to the \textit{synthesis} point of view when employing sparsity-based methods. That is, it corresponds to a generative model, where the signals are assumed to be well-representable by a sparse combination of the dictionary elements. The complementary concept of employing sparsity-based methods relying on the \textit{analysis} point of view, where the considered signals are assumed to be sparse \textit{after} the application of a sparsifying transform exists as well and has also been considered for MR image reconstruction, both on a patch-based level \cite{seibert2014separable, ravishankar2015efficient} as well as using convolutional transforms \cite{kofler2022cdl}.

\subsection{Data-Driven Dictionary-Based  Regularization Methods for Quantitative MRI}\label{sec:data_driven_qmri}

Dictionary learning can be employed for quantitative MRI reconstruction as well. For example, one can follow a  decoupling approach and use the problem formulation \eqref{eq:dlmri_psi_ada_learned} in \eqref{eq:qmri_variational_problem_splitted2} as done in \cite{doneva2010compressed}, i.e.\
\begin{equation}\label{eq:qdlmri_variational_problem_splitted} 
    \left \{
\begin{aligned}
&\underset{\PP}{\min}\;  \frac{1}{2}\| q_{\mathcal{M}}(\PP) - \XX \|_2^2, \\ 
&\text{where}\;\; \XX = \underset{\XX^\prime, \mathbf{\Psi} \in \mathcal{D}_{d,K},\{\boldsymbol{\gamma}_j\}_j}{\argmin} \frac{1}{2}\| \Au  \XX^\prime -\YY\|_2^2 + &\frac{\lambda}{2} \sum_{j=1}^M \| \Rd_j \XX^\prime - \mathbf{\Psi} \boldsymbol{\gamma}_j \|_2^2 \\ & &\text{s.t.} \quad  \forall j: \| \boldsymbol{\gamma}_j\|_0 \leq S.
\end{aligned}\right.
\end{equation}
Alternatively, a coupling approach can be adopted by appropriately adapting the problem formulation \eqref{eq:dlmri_psi_ada_learned} to ensure that the patches of the quantitative parameters are sparse with respect to a learned dictionary, i.e.\ 
\begin{eqnarray}\label{eq:qdlmri_psi_ada_learned}
\underset{\PP, \mathbf{\Psi} \in \mathcal{D}_{d,K}, \{\boldsymbol{\gamma}_j\}_j}{\mathrm{min}} \frac{1}{2}\| (\Au \circ q_{\mathcal{M}})(\PP) -\YY\|_2^2 + \frac{\lambda}{2} \sum_{j=1}^M \| \Rd_j \PP - \mathbf{\Psi} \boldsymbol{\gamma}_j \|_2^2 \nonumber \\ \quad \text{such that} \quad   \forall j: \| \boldsymbol{\gamma}_j\|_0 \leq S,
\end{eqnarray}
see \cite{kofler2023quantitative}. Thereby, we note that from a computational point of view, problem \eqref{eq:qdlmri_psi_ada_learned} is much more attractive than \eqref{eq:qdlmri_variational_problem_splitted}, since in general, $N_p\ll Q$ and therefore, there are only \textit{few} quantitative parameters $\PP$ for which one needs to repeatedly apply a sparse coding algorithm in contrast to \eqref{eq:dlmri_psi_ada_learned}, where the sparse coding has to be repeatedly performed for the \text{many} qualitative images $\XX$. This advantage however comes at the cost of the non-linearity of the entire forward operator $\Au \circ q_{\mathcal{M}}$ in \eqref{eq:qdlmri_psi_ada_learned} in contrast to the linear one $\Au$ in the constraint in \eqref{eq:qdlmri_variational_problem_splitted}. This must be appropriately addressed, either using non-linear gradient descent methods or variable splitting strategies as in \cite{kofler2023quantitative}.

\section{Neural Networks for MR Image Reconstruction: A Survey of Different Approaches}\label{sec:nns_for_mri}

Although the previously discussed methods are data-driven in the sense that the sparsifying transforms are learned from data, the underlying principle that these methods use for regularization, i.e.\ the concept of sparsity, is still hand-crafted. To overcome this limitation, one can adopt more powerful models that do not rely on specific hand-crafted concepts but are rather designed to perform a certain task, e.g.\ denoising or the reduction of artifacts.

In addition, from the formulations \eqref{eq:dico_learning_problem_L0} and \eqref{eq:cd_learning_problem_L1}, we see an inherent conceptual limitation of these learning paradigms: The problem formulations for learning the dictionary ignore the underlying physics that is responsible for the generation of signals/images, i.e.\ the forward model never appears in the problem formulation for learning. Thus, with these types of approaches, the rationale behind the learning is often to \textit{first learn} the regularization method, and only \textit{then} to reconstruct images/quantitative parameters for a specific inverse problem of interest. This particular aspect can be overcome by employing Neural Networks (NNs).

NNs are highly versatile functions that can be employed in many different ways for the regularization of inverse problems. In particular, so-called model-based (also sometimes called model-aware or physics-informed) NNs form methods for which the training conceptually 
differs from the previous we have seen for example in \eqref{eq:dico_learning_problem_L0} and \eqref{eq:cd_learning_problem_L1}. Instead of decoupling the learning process from the reconstruction method, model-aware NNs offer the possibility to be \textit{learned to reconstruct} and for that reason, they nowadays define the state-of-the-art in image reconstruction across all imaging modalities, not only for MRI.

We mention that naturally many reviews exist in the topic of NNs-based regularization methods for inverse problems including MRI \cite{arridge_solving_2019, kamilov2023plug, monga2021algorithm, mukherjee_learned_2023, shlezinger2023model}, see also \cite{habring2023neural} and the references within. However, a comprehensive review and categorization of how neural networks can produce state-of-the-art results in quantitative MRI is missing from the literature. Here we aim to cover this gap.

Before discussing model-based NNs in more detail, we briefly introduce the basic concepts of deep learning that are necessary to emphasize the differences between model-aware learning and model-agnostic learning of NNs.

\subsection{Neural Networks in a Nutshell}

\subsubsection{Feed-Forward Neural Networks}

For the sake of simplicity, we restrict our exposition in this section to so-called feed-forward NNs. More advanced network architectures are also used in practice, see also the next Section \ref{sec:Unet}.
In its most general form, an NN is simply a function between two finite-dimensional Euclidean spaces  $V, W$, parameterized by a set of  parameters $\Theta\in\R^{\ell}$ that are to be tuned:
\begin{eqnarray}\label{eq:nn}
    f_{\Theta}: V &\longrightarrow &W, \\
                v &\longmapsto &w:=f_{\Theta}(v).
\end{eqnarray}
Typically, $f_{\Theta}$ is constructed as a composition of multiple functions $f_{\Theta_i}$, $i=1,\ldots, N_{\mathrm{layers}}$ - the so-called layers of the NN - which comprise relatively simple functions (affine). These affine maps are represented by matrices and vectors denoted by $\mathbf{W}_i\in\R^{n_{i}\times n_{i-1}}$ and $\mathbf{b}_i\in \R^{n_{i}}$  respectively, often referred to as the \textit{weights} and \textit{biases} of the network. 
These are each followed by non-linearities $\sigma_i$, the so-called activation functions. Thus, the functions $f_{\Theta_i}$ are typically of the form $f_{\Theta_i}(\,\cdot\,):= \sigma_i( \mathbf{W}_i \cdot + \mathbf{b}_i)$, 
where $\Theta_i:=\{\mathbf{W}_i,\mathbf{b}_i\}$ and the application of $\sigma_i$ is understood to be  component-wise.  We note that the activation functions can also be parameterized by some learnable parameters as well. A network $f_{\Theta}$ with $N_{\mathrm{layers}}$  is formed via a successive composition of the functions $f_{\Theta_i}$, $i=1,\ldots,N_{\mathrm{layers}}$:
\begin{equation}\label{eq:nn_layers}
    f_{\Theta} = \bigcirc_{i=1}^{N_{\mathrm{layers}}} f_{\Theta_i}, \text{ where }\, \Theta:=\bigcup_{i=1}^{N_{\mathrm{layers}}} \Theta_i.
\end{equation}

In general, there exist a large variety of different network architectures (i.e.\ number of layers, size and structure of weight matrices, choice of activation functions) and we will abstain from discussing them in detail. Here, we only mention that the structure of the weights $\mathbf{W}_i$ and the biases $\mathbf{b}_i$ is often determined  by the considered task as well as by the available computational resources. For example, in the context of image processing tasks as image segmentation or image denoising, the weight matrices $\mathbf{W}_i$ are often chosen such that they correspond to convolutions. For image classification tasks, in the last layers of the NNs, which are used to perform a classification of the features extracted by the previous convolutional layers, the matrices $\mathbf{W}_i$ usually correspond to full-matrices. 
For a detailed introduction to NNs and to modern architectures, we refer to \cite{goodfellow2016deep}.

The popularity of NNs stems from the fact that they are able to approximate arbitrary well (uniformly in compact sets) any continuous function (Universal Approximation Theorem, \cite{hornik1989multilayer}); 
a property that is extended to more general functions under an appropriate topology for the approximation \cite{pinkus_1999}. 
NNs are trained to have some desirable effect on a set of data (training set) and end up having a similar effect on unseen data as well (test set). This effect is typically quantified with the help of a loss-function that is minimized during training, which we describe in the next subsection. This means that for practical applications, NNs can be trained such that they have remarkable generalization capabilities beyond their training data.

\subsubsection{The U-Net}\label{sec:Unet}

Although the focus of this chapter is not the in-depth discussion of state-of-the-art network architectures but rather their application within the task of MR reconstruction problems, we briefly describe the U-Net \cite{ronneberger2015u} as it is without doubt one of the most employed architectures.

\begin{figure}[t]
    \centering
    \includegraphics[width=0.7\linewidth]{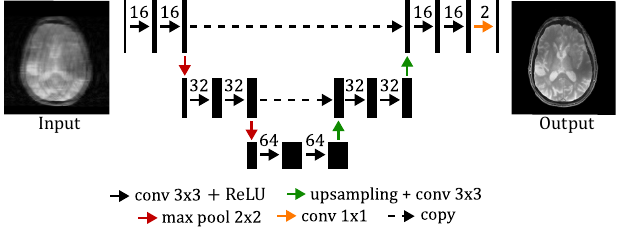}
    \caption{An example of a U-Net \cite{ronneberger2015u} for reducing artifacts and noise from the zero-filled reconstruction of undersampled $k$-space data. This particular example of a U-Net consists of three resolution scales. The numbers over the convolution arrows signify varying numbers of convolution filters used. Here, each convolutional block consists of two convolutional layers followed by a rectified linear unit (ReLU) activation function. The last layer maps to the real and imaginary part of the complex-valued output image.}
    \label{fig:unet}
\end{figure}

The U-Net (see Figure  \ref{fig:unet}) can be regarded as a multi-scale convolutional network, comprised of blocks of convolution layers (with learnable weights) and accompanying activation functions, spatial downsampling, and upsampling steps. In this context, downsampling can be achieved, for instance, by taking the maximum value over a window (\textit{max-pooling}), while upsampling can be implemented as bilinear interpolation. The network incorporates \textit{skip-connections} as shortcuts for information at different resolution scales. Commonly, the sequence of convolutional blocks and downsampling operations before these skip connections is regarded as the \textit{encoder}, and the upsampling operations and convolutional blocks after the skip connections as the \textit{decoder}. Based on this skeleton, many extensions and variations have been proposed \cite{ye2018deep,zhou2018unet++,oktay2018attention,qin2020u2}.

\subsection{Training of Neural Networks}
Training an NN refers to optimizing its trainable parameters $\Theta$ on a \textit{training set}. For instance, in the case of {\it supervised learning} this training set consists of input-target pairs $\mathcal{D}:=\{(v_i,w_i)_{i=1}^{N_{\mathrm{train}}}\}$. Training is in general achieved by minimizing a so-called loss-function using (stochastic) gradient descent-type schemes. There, one constructs a sequence of parameters $\{\Theta_k\}_k$ by
\begin{equation}\label{eq:grad_descent}
    \Theta_{k+1}:=\Theta_k - \tau_k (\nabla_{\Theta} \mathcal{L})(\Theta_{k}),
\end{equation}
where $\mathcal{L}$ denotes a loss-function of the form
\begin{equation}\label{eq:loss_fct}
    \mathcal{L}(\Theta):=\frac{1}{|\mathcal{D}|}\sum_{(v,w)\in\mathcal{D}} l(f_{\Theta}(v), w) + r(\Theta).
\end{equation}
Here $\{\tau_k\}_{k}$ is a sequence of positive step-size parameters also known as the learning rate, which can either be chosen a priorr, or also adaptively varied, see e.g.\ \cite{takase2018effective}. The function $l(\,\cdot\,,\,\cdot\,)$ is a discrepancy metric that is appropriately chosen with respect to the considered task, e.g.\ the mean squared error (MSE) for image processing tasks, the cross-entropy for classification or the dice similarity-based loss for segmentation tasks.
The term $r(\cdot)$ is a regularization term that can be used to avoid overfitting of the network which we explain in the next paragraph.

For model development, one usually divides all available data into three portions: the \textit{training, validation} and \textit{test} data. The training data is the dataset used to minimize \eqref{eq:loss_fct}, while on the validation dataset, one monitors the values of the loss-function $l(\,\cdot\,,\,\cdot\,)$ during training. This step is needed to avoid overfitting, which refers to the phenomenon where the loss-function evaluated on the training set is much smaller than its corresponding evaluation on the validation set, i.e. the network does not generalize well beyond the training data. This often occurs - and hence the role of $r(\cdot)$ becomes important - when $|\Theta|$ is large (e.g.\ in the range of millions) while $|\mathcal{D}|$ is relatively small (e.g.\ in the range of hundreds). The validation set is additionally also used to test different model configurations, e.g.\ different number of layers, learning rates, etc. Finally, the model's parameters for which the loss-function evaluated over the entire validation set is the smallest, is applied to the test set. Note that the meaning of the terms validation and test set is sometimes interchanged in the literature.  

Instead of only relying on plain gradient descent schemes as described in \eqref{eq:grad_descent}, much work has been put into the development of non-linear optimization routines that work well in the high-parameter setting of NNs.
Well-known optimization algorithms are RMSprop \cite{hinton2012neural}, AdaDelta \cite{zeiler2012adadelta} or Adam \cite{kingma2014adam}, to only name a few, see \cite{ruder2016overview} for an overview and comparison.
Additionally, for computational reasons, one typically rarely utilizes the entire dataset $\mathcal{D}$ to compute the gradient of $\mathcal{L}$, but employs stochastic gradient schemes, i.e.\ during training, the sum in \eqref{eq:loss_fct} is approximated by a (typically) small number of samples randomly chosen from $\mathcal{D}$ \cite{nocedal1999numerical}, also referred to as mini-batches. This strategy not only is attractive from a computational point of view but can also be seen as a form of implicit regularization, see \cite{wilson2003general, smith2021on}.

There exist different ways of optimizing the NNs parameters, depending on the available training data as well as on the specific task. Thus, we next discuss different options for specifically training NNs for MR image reconstruction in more detail.

\subsubsection{Supervised Training}
For supervised training, one requires access to a set of input-target pairs and the goal is to train the network's weights such that the network can be used to predict new estimates of the target quantities.

Recall the problem formulations in \eqref{eq:mri_forward_model} and \eqref{eq:qmri_forward_model}. To be able to apply supervised training for qualitative and quantitative MR reconstruction networks, let us assume to have the following datasets available:
\begin{equation}\label{eq:sup_mri_training_data}
    \mathcal{D}_{\mathrm{sup}}^{\mathrm{MRI}}:= \{ (\YY^i, \XX_{\mathrm{true}}^i)_{i=1}^{N_{\mathrm{train}}}\,\,  \big|  \,\,  \,\YY^i:= \Au \XX_{\mathrm{true}}^i + \mathbf{e}^i, \; i=1,\ldots,N_{\mathrm{train}}   \},
\end{equation}
and
\begin{align}
    \mathcal{D}_{\mathrm{sup}}^{\mathrm{qMRI}_1}&:= \{ (\YY^i, \PP_{\mathrm{true}}^i)_{i=1}^{N_{\mathrm{train}}}\,\,  \big|  \,\,  \,\YY^i:= \big(\Au \circ q_{\mathcal{M}}\big) (\PP_{\mathrm{true}}^i) + \mathbf{e}^i, \; i=1,\ldots,N_{\mathrm{train}}   \}\label{eq:sup_qmri_training_data1},\\
    \mathcal{D}_{\mathrm{sup}}^{\mathrm{qMRI}_2}&:= \{ (\XX^i, \PP_{\mathrm{true}}^i)_{i=1}^{N_{\mathrm{train}}}\,\,  \big|  \,\,  \,\XX^i:= q_{\mathcal{M}}(\PP_{\mathrm{true}}^i) + \boldsymbol{\eta}^i, \; i=1,\ldots,N_{\mathrm{train}}   \}.\label{eq:sup_qmri_training_data2}
\end{align}
Then, for the respective MR and quantitative MR reconstruction problems, training an appropriate network $f_{\Theta}$ in a \textit{supervised} manner corresponds to minimizing the following loss-functions:
\begin{equation}\label{eq:mri_sup_training}
    \mathcal{L}_{\mathrm{sup}}^{\mathrm{MRI}}(\Theta):=\frac{1}{|\mathcal{D}_{\mathrm{sup}}^{\mathrm{MRI}}|}\sum_{(\YY,\, \XX_{\mathrm{true}})\in\mathcal{D}_{\mathrm{sup}}^{\mathrm{MRI}}} l\big(f_{\Theta}(\YY),   \XX_{\mathrm{true}}\big) + r(\Theta),
\end{equation}
and
\begin{equation}\label{eq:qmri_sup_training1}
    \mathcal{L}_{\mathrm{sup}}^{\mathrm{qMRI}_1}(\Theta):=\frac{1}{|\mathcal{D}_{\mathrm{sup}}^{\mathrm{qMRI}_1}|}\sum_{(\YY,\,\PP_{\mathrm{true}})\in\mathcal{D}_{\mathrm{sup}}^{\mathrm{qMRI}_1}} \big(f_{\Theta}(\YY),   \PP_{\mathrm{true}}\big) + r(\Theta),
\end{equation}
or
\begin{equation}\label{eq:qmri_sup_training2}
    \mathcal{L}_{\mathrm{sup}}^{\mathrm{qMRI}_2}(\Theta):=\frac{1}{|\mathcal{D}_{\mathrm{sup}}^{\mathrm{qMRI}_2}|}\sum_{(\XX,\,\PP_{\mathrm{true}})\in\mathcal{D}_{\mathrm{sup}}^{\mathrm{qMRI}_2}} \big(f_{\Theta}(\XX),   \PP_{\mathrm{true}}\big) + r(\Theta).
\end{equation}
Note that in \eqref{eq:qmri_sup_training1}, $f_{\Theta}$ is assumed to be a learned network that maps $k$-space data directly to quantitative parameters, while in \eqref{eq:qmri_sup_training2}, $f_{\Theta}$ maps the series of qualitative images to the quantitative parameters. Different networks in these forms will be discussed later in Section \ref{sec:nns_for_mri}.

\subsubsection{Self-Supervised Training}

In many applications, obtaining target images or quantitative parameters to be used as ground-truth target data can be a challenging or even impossible task. For example, when imaging dynamic processes in cardiac MRI, there inherently exists a trade-off between temporal and spatial resolution, and consequently, the $k$-space always needs to be undersampled. The same holds for quantitative imaging, where the measurement process is repeated multiple times to obtain the $k$-space of the different $Q$ qualitative images from which the quantitative parameters can then be estimated.

In self-supervised learning, one exploits the structure of the type of problem we are considering. Recall that we are dealing with data generated according to the forward model in  \eqref{eq:mri_forward_model}. This means that, if not provided with a dataset of the type \eqref{eq:sup_mri_training_data}, \eqref{eq:sup_qmri_training_data1} or \eqref{eq:sup_qmri_training_data2}, i.e.\ without access to target images or target quantitative parameters, we can make use of the known forward model and apply it to the quantity predicted by our reconstruction method $f_{\Theta}$. For example, for qualitative MR image reconstruction, an example of an adaptable loss-function could be given by
\begin{equation}\label{eq:loss_fct_self_sup_mri}
    \mathcal{L}_{\text{self-sup}}(\Theta):=\frac{1}{|\mathcal{D}|}\sum_{i=1}^{N_{\mathrm{train}}} l\bigg( \Au f_{\Theta}(\YY^i) \, , \, \YY^i\bigg) + r(\Theta),
\end{equation}
while for quantitative MR image reconstruction, there are two different possibilities, since the entire forward model $\Au \circ q_{\mathcal{M}}$ is a composition of two operators. For the case where target qualitative images are available, one can use 
\begin{equation}\label{eq:loss_fct_self_sup_qmri1}
    \mathcal{L}_{\text{self-sup}}^{\text{image}}(\Theta):=\frac{1}{|\mathcal{D}|}\sum_{i=1}^{N_{\mathrm{train}}} l\bigg( q_{\mathcal{M}}\big(f_{\Theta}(\YY^i)\big) \, , \, \XX_{\mathrm{true}}^i \bigg) + r(\Theta),
\end{equation}
while if only the measured $k$-space data is available, one can minimize 
\begin{equation}\label{eq:loss_fct_self_sup_qmri2}
    \mathcal{L}_{\text{self-sup}}^{\text{k-space}}(\Theta):=\frac{1}{|\mathcal{D}|}\sum_{i=1}^{N_{\mathrm{train}}} l\Big( \big(\Au \circ q_{\mathcal{M}}\big)\big(f_{\Theta}(\YY^i)\big) \, , \, \YY^i \Big) + r(\Theta).
\end{equation}
Note that,  for the choice of $l(\, \cdot \, ,\, \cdot \,) = \|\, \cdot \, - \, \cdot \, \|_2^2$ for instance, we have 
\begin{align}
\|\Au f_{\Theta}(\YY)-\YY\|_{2}
&\leq\| \Au f_{\Theta}(\YY)-\Au \XX_{\mathrm{true}} - \EE \|_2\nonumber\\ 
&\leq \| \Au \|\cdot \| f_{\Theta}(\YY) - \XX_{\mathrm{true}} \|_2 + \|\EE\|_{2}.\label{eq:ss_vs_sup_ineq_mri}
\end{align}
The inequality \eqref{eq:ss_vs_sup_ineq_mri} suggests that if a network has been trained in a supervised learning manner such that $f_{\Theta}(\YY)\approx \XX_{\mathrm{true}}$, i.e.\  the loss-function in \eqref{eq:mri_sup_training} is relatively small, then,  provided that the amount of noise is small, the same network will also make the self-supervised loss-function \eqref{eq:loss_fct_self_sup_mri} small. Thus, it is expected that  $f_{\Theta}$ will work roughly at least as well as any network that has been trained using \eqref{eq:loss_fct_self_sup_mri}. Note that an argument in the opposite direction would require a reverse inequality of \eqref{eq:ss_vs_sup_ineq_mri} which would involve $\|\Au^{-1}\|$. Even in the case where this quantity is well-defined, it would be large due to ill-posedness of the inverse problem. Hence, one expects that supervised learning approaches produce in general superior results to self-supervised ones.

Similarly, for quantitative MR, by exploiting the fact that the non-linear signal model $q_{\mathcal{M}}$ is typically Lipschitz-continuous \cite{Dong_2019}, and by denoting the Lipschitz-constant by $L_{q_{\mathcal{M}}}$  the following inequality holds in view of \eqref{coupled_ip_1} and \eqref{coupled_ip_2}
\begin{align}\label{eq:ss_vs_sup_ineq_qmri}
\|\big(\Au \circ q_{\mathcal{M}}\big)\big(f_{\Theta}(\YY)\big)-\YY\|_{2}
&\le\| \big(\Au \circ q_{\mathcal{M}}\big)\big(f_{\Theta}(\YY)\big)- \big(\Au \circ q_{\mathcal{M}}\big) (\PP_{\mathrm{true}}) - \Au\boldsymbol{\eta}-\EE\|_2 \nonumber\\ 
&\leq  \| \Au \|\cdot \| q_{\mathcal{M}} (\PP_{\mathrm{true}}) - q_{\mathcal{M}}\big(f_{\Theta}(\YY)\big) \|_2  + \|\Au\boldsymbol{\eta}\|_{2}+\|\EE\|_{2} \nonumber \\
&\leq\| \Au \| \cdot L_{q_{\mathcal{M}}} \cdot \| \PP_{\mathrm{true}} - f_{\Theta}(\YY) \|_2 + \|\Au\boldsymbol{\eta}\|_{2}+\|\EE\|_{2},
\end{align}
which shows the analogous for learned reconstruction methods for quantitative MRI.

\subsubsection{Self-Supervised Training by Data-Undersampling}

Another theoretically founded possibility is to consider self-supervised training by data-undersampling (SSDU) \cite{yaman2020self} which we briefly discuss here.
SSDU was recently shown to be closely connected to a concept developed in the image denoising community \cite{millard2023theoretical}, i.e.\  the Noisier2Noise \cite{moran2020noisier2noise} framework, which extends the Noise2Noise \cite{lehtinen2018noise2noise} framework.

In the Noise2Noise framework \cite{lehtinen2018noise2noise}, it was shown that denoising NNs can also be trained without access to paired input-target image pairs. Instead, it is sufficient to have access to input pairs that are corrupted by the statistically same type of noise. While this approach does not require clean target images, it still requires access to the same input image that is corrupted by a different instance of the statistically same noise. This approach is interesting for the case where one can afford to easily repeat a measurement, which is clearly not the case for MRI. The Noise2Self \cite{batson2019noise2self} and Noisier2Noise \cite{moran2020noisier2noise} frameworks can be adopted to overcome this limitation. 

Noise2Self allows to perform image denoising by self-supervision, i.e.\ the same given noisy data can be used as a target during training. This is achieved by defining two disjoint partitions of the pixels of the image, utilizing one portion of the pixel as the input to the network, and by computing the loss between the prediction and the given noisy image when restricted to the second partition of the pixels.
In Noisier2Noise, it is shown that it is sufficient to further corrupt a given noisy image and to map this noisier sample to its original noisy one - hence the name Noisier2Noise. Then, at inference time, it is possible to employ the trained network to obtain an estimate of the unknown ground-truth image after a linear correction. 

SSDU \cite{yaman2020self} builds on these ideas by the following. Let $\Au$ denote a forward operator of the form in \eqref{eq:mri_forward_model} with $k$-space coefficients sampled at the frequencies defined by the set $I$. In SSDU, one now chooses a split of the coefficients $I:=I_1 \cup I_2$ with $I_1 \cap I_2 = \emptyset$ and uses the portion $\YY_{I_{1}}$ of the $k$-space data to reconstruct an estimate of the image and the remaining $\YY_{I_{2}}$ for the computation of a loss-function. The loss-function minimized when applying SSDU is of the form
\begin{equation}\label{eq:loss_fct_ssdu_mri}
    \mathcal{L}_{\text{SSDU}}^{\text{MRI}}(\Theta):=\frac{1}{|\mathcal{D}|}\sum_{i=1}^{N_{\mathrm{train}}} l\big( \Ad_{I_2} f_{\Theta}(\YY_{I_1}^i)\, , \, \YY_{I_2}^i\big) + r(\Theta).
\end{equation}
In \cite{millard2023theoretical}, the theoretical connection between SSDU and the Noisier2Noise framework is described in more detail. More precisely, it also gives insight into how the partition of $I$ should take place to fulfill the requirements of Noisier2Noise. Additionally, \cite{millard2023theoretical} suggests a particular choice of loss-function $l(\,\cdot\, , \,\cdot\, )$ that uses a weighting to compensate for the sampling and partitioning densities, see the respective paper for details.

Similarly,  for a quantitative MR problem of the form \eqref{coupled_ip_2} with only access to qualitative images, one can for example randomly split an image $\XX$ into two images with disjoint support, i.e.\ $\XX = \XX_{I_1} + \XX_{I_2}$, where, by slight abuse of notation, we denote $\XX_J[j]=\XX[j]$ if $j\in J$ and $\XX_J[j]=0$ if $j\notin J$. Then, a Noise2Self-inspired loss-function
\cite{zimmermann2023ismrm}  for a quantitative MR reconstruction network $f_{\Theta}$ that maps qualitative images to quantitative parameters can be constructed as 
\begin{equation}\label{eq:loss_fct_ssdu_qmri1}
    \mathcal{L}_{\text{SSDU}}^{\text{qMRI}_1}(\Theta):=\frac{1}{|\mathcal{D}|}\sum_{i=1}^{N_{\mathrm{train}}} l\bigg(  \big(q_{\mathcal{M}}\,f_{\Theta}(\XX_{I_1}^i)\big)_{I_2}\, , \, \XX_{I_2}^i\bigg) + r(\Theta).
\end{equation}
Last, if only the $k$-space data is available, one can use
\begin{equation}\label{eq:loss_fct_ssdu_qmri2}
    \mathcal{L}_{\text{SSDU}}^{\text{qMRI}_2}(\Theta):=\frac{1}{|\mathcal{D}|}\sum_{i=1}^{N_{\mathrm{train}}} l\bigg(  \big(\Ad_{I_2} \circ q_{\mathcal{M}}\big) \big(f_{\Theta}(\YY_{I_1}^i)\big)\, , \, \YY_{I_2}^i\bigg) + r(\Theta)
\end{equation}
and employ the SSDU-technique as described above.

\subsection{Publicly Available Datasets}

Regardless of the loss-functions employed for training, deep learning-based approaches typically require the availability of relatively large datasets to be properly trained. However, obtaining large, high-quality datasets that are adaptable for algorithmic development can be a challenging task. Additionally, the data used for training can often not be shared among institutions due to data-privacy related issues, which hinders the comparability of deep learning-based approaches. 

To overcome these limitations, in the last years, several initiatives have taken place aimed at collecting, organizing, and publishing large amounts of high-quality datasets to be used for training and to facilitate the fair comparison of data-driven reconstruction methods.

Additionally, in the last few years, several MRI reconstruction challenges have been organized and attracted the attention of many researchers around the globe. 

Typically, the challenges are organized in a similar way. In the training phase of the challenges,  research teams register for the challenge and obtain access to (typically paired) input-target images, where the goal is to develop methods that are subsequently applied to the test data, for which the target images are usually withheld by the challenge organizers.

Then, the challenge organizers report pre-defined image metrics that the participating teams obtained on the test set. The metrics typically consist of error-based metrics such as the normalized root mean squared error (nRMSE) or the peak signal-to-noise ratio (PSNR) and structural similarity index measure (SSIM) \cite{wang2004image}.

Here, we only briefly the availability of a few of the most famous datasets, e.g.\ the fastMRI dataset \cite{zbontar2018fastmri} and the associated fastMRI reconstruction challenge \cite{muckley2021results} (for brain and knee MRI), the Calgary-Campinas initiative \cite{souza2018open} and the associated multi-coil MRI reconstruction challenge \cite{beauferris2022multi} (for brain MRI), the CMRxRecon dataset and its associated challenge \cite{wang2023cmrxrecon} (for cardiac cine MRI) and the K2S challenge for joint knee MR image reconstruction and segmentation \cite{tolpadi2023k2s}.

\subsection{Categorization of Neural Networks Methods for MRI}

Even though the introduction of neural networks in image reconstruction led to the production of state-of-the-art, concerns have been raised about their ``black box'' nature especially when they are used in a direct model-inversion fashion, see the next subsection. For the latter, there is little interpretability of the reconstruction results, and potential instabilities and unexpected artifacts could interfere with the diagnostic procedure \cite{antun2020instabilities}. This has led to the development of families of methods that combine neural networks with techniques from the traditional hand-crafted approaches mentioned in the previous sections. These approaches integrate the interpretability and convergence guarantees of classical variational methods with the flexibility and adaptability of neural networks. In what follows, we will be mainly focusing on this family of methods for MRI and qMRI reconstruction. 
Additionally, we will further categorize these combined methods based on whether they are model-aware or model-agnostic during the training and the inference phase according to Table \ref{tbl:four_categories}.

\begin{table}[ht]
\renewcommand{\arraystretch}{1.3}
\noindent
\setlength\tabcolsep{3pt}
\begin{tabular}[t]{|m{0.13\textwidth} | m{0.40\textwidth} | m{0.40\textwidth} | }
\hhline{~--}
\multicolumn{1}{c|}{} 
&
\multicolumn{1}{c|}{\cellcolor{gray!30}\textbf{Model-agnostic}} 
&
\multicolumn{1}{c|}{\cellcolor{gray!30}\textbf{Model-aware}}\\
\hline
\cellcolor{gray!30}
\textbf{Training phase} 
& 
\footnotesize{The {\it forward model  is not involved} in the training of the neural network apart perhaps in its first layer, see Section \ref{sec:model_agnostic_learning_mri}.}
& \footnotesize{The {\it forward model is involved} in the training of the neural network in a non-trivial way, i.e. not only to simulate data.}\\
\hline
\cellcolor{gray!30}
\textbf{Inference phase} 
& \footnotesize{After (any) training is completed, {\it the method does not use the forward model to reconstruct} the image/quantitative maps from given data.}
&  \footnotesize{After (any) training is completed, {\it the method uses the forward model to reconstruct} the image/quantitative maps from given data, see Section \ref{sec:model_aware}.}\\
\hline
\end{tabular}
\caption{Differentiation of data-driven reconstruction methods based on the used learning and reconstruction paradigm.}
\label{tbl:four_categories}
\end{table}

{\it Remark}: A clarification over the term ``non-trivial way'' in the description of the model-aware methods in Table \ref{tbl:four_categories} is in order. There is a wide spectrum of ways that the forward model can be involved in the training of the neural network resulting in different degrees of interpretability. Here, interpretability refers among others to the degree of similarity of the network's architecture with classical model-based algorithmic schemes. For instance, on one side of the spectrum we have post-processing methods, see Section \ref{sec:post_processing} below, where $\Au$ (more precisely its Moore-Penrose pseudoinverse $\Au^\dagger$) appears only as a very first component of the network. Thus, knowledge about the model is only used to transfer the data to the image space but does not result in any further interpretability or data-consistency. On the other hand, the interpretability of a method could be for instance strengthened if the NN contains repeated components (e.g.\ layers) representing, for example, gradient descent steps of the data-fidelity term, i.e. $\XX_{k+1}=\XX_{k}-\tau \Ad_{I}^{\herm}(\Ad_{I}\XX_{k}-\YY)$. In fact, replacing components of iterative reconstruction schemes with data-driven blocks forms a convenient way to derive meaningful model-aware methods, see Section \ref{sec:model_aware}.

\subsection{Neural Networks for Qualitative MRI}\label{seq:nn_for_mri}

We first start again with the exposition of NNs-based methods to be used as regularization methods for qualitative MR image reconstruction problems, since, as we have seen before in \eqref{eq:qmri_variational_problem_splitted1} and \eqref{eq:qmri_variational_problem_splitted2}, they can be used for reconstructing qualitative MR images prior to a non-linear fit to obtain the quantitative parameters. We begin with NNs-based methods for regularization that are model-agnostic in the training phase, i.e.\ they fall in the category of the left column of Table \ref{tbl:four_categories}. 

\subsubsection{Direct Model-Inversion}

\begin{figure}[h!]
    \centering
    \includegraphics[width=0.9\linewidth]{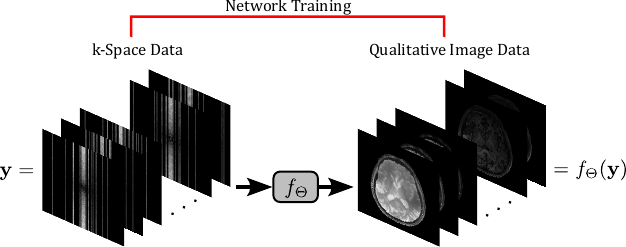}
    \caption{
    An example of a model-agnostic NN given by a learned direct inversion of the forward model $\Au$ is shown for a brain MRI example. For qualitative MR image reconstruction, the method named AUTOMAP (AUtomated TransfOrm by Manifold APproximation) can be found in \cite{zhu2018image}.
    }  
    \label{fig:full_inversion}
\end{figure}

Here, we start with the perhaps simplest way to employ NNs for the task of image reconstruction.
The simplest approach is to completely ignore the knowledge about the structure of the forward problem in \eqref{eq:mri_forward_model} and to learn a direct inversion from the observed data to an estimate of the target image. The procedure is visually illustrated in Figure \ref{fig:full_inversion}.
An early idea for this approach was originally presented in \cite{paschalis2004tomographic} and \cite{argyrou2012tomographic} for some relatively small-scale problems, and later on re-discovered for larger problems employing the use of large GPUs for training \cite{zhu2018image}.

From Figure \ref{fig:full_inversion}, one can see why this approach is both model-agnostic in the learning phase as well as at inference time and therefore, falls in the left column of Table \ref{tbl:four_categories}. Clearly, the learning phase is model-agnostic as the NN $f_{\Theta}$ entirely ignores the physical model $\Au$ that is responsible for the generation of the measured $k$-space data. No information about the MRI device or the data-acquisistion process is utilized for the design of the network $f_{\Theta}$.

In addition, the approach presented in \cite{zhu2018image} has some limitations regarding its applicability. First of all, because the first operations of  $f_{\Theta}$  consist of several \textit{fully-connected} layers that aim at learning the inversion of the Fourier-transform, its applicability is restricted to images of the same dimensions that $f_{\Theta}$ was trained with, in contrast to convolutional layers. Additionally, the number of parameters to be trained is extremely large by construction, and the fact that the network needs to entirely learn a mapping from two different spaces makes it challenging to apply the method to large-scale problems.

\subsubsection{Post-Processing Methods}\label{sec:post_processing}

Instead of learning the entire inversion of the forward model, it might be beneficial to learn to map an initial reconstruction that can be trivially obtained from the measured data to an estimate of its unknown ground-truth image.

From the observed $k$-space data $\YY$, an initial guess of the qualitative MR image $\XX_{\mathrm{true}}$ can for example be obtained by simply applying the adjoint $\Au^\herm$. In the case where $\Au$ is injective, e.g.\ for the case of MR imaging with multiple receiver coils, one can approximately solve the normal equations
\begin{equation}\label{eq:normal_eqs}
    \Au^\herm \Au \XX = \Au^\herm \YY,
\end{equation}
and obtain $\XX^\dagger:=\Au^\dagger \YY$, where $\Au^\dagger:= (\Au^\herm \Au)^{-1} \Au^\herm$ denotes the Moore-Penrose pseudoinverse of $\Au$. Thereby, the solution is typically obtained by solving a system of linear equations with an iterative method, e.g.\ the conjugate gradient (CG) method, rather than inverting the operator $\Au^\herm \Au$. 
Additionally, as for non-Cartesian acquisitions, pre-conditioning in $k$-space is often applied \cite{pruessmann2001advances}, one can as well use the density-compensated $k$-space data to obtain an initial solution, or by choosing a weighted $\ell_2$-norm as data-discrepancy measure, and solve the pre-conditioned normal equations
\begin{equation}\label{eq:normal_eqs_dcomp}
    \Au^\herm\, \Wd\, \Au \XX = \Au^\herm\, \Wd\, \YY.
\end{equation}

\begin{figure}
    \centering
    \includegraphics[width=\linewidth]{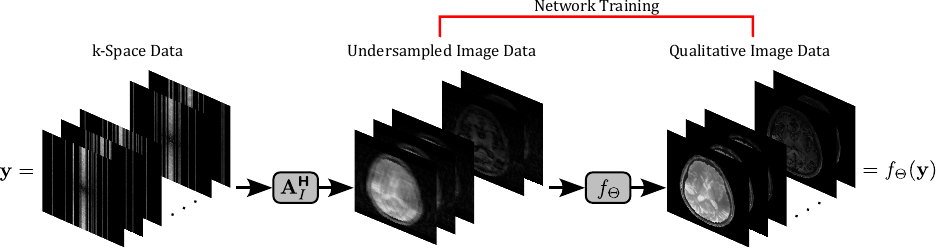}
    \caption{
    An example of a model-agnostic NN given by a learned post-processing method that reduces noise and artifacts that are present in the initial zero-filled reconstruction $\XX_0:=\Au^\herm \YY$.}
    \label{fig:post_processing}
\end{figure}

Depending on the exact set-up of the considered inverse problem, the initial solution contains strong artifacts and/or noise. 
Figure \ref{fig:post_processing} illustrates again the process for the post-processing of a zero-filled reconstruction with $f_{\Theta}$ for a brain MR example.

Similar to Figure \ref{fig:full_inversion}, we see that the training of the network $f_{\Theta}$ is model-agnostic as well because neither the forward nor the adjoint model are used to construct the network.  The operation of $\Au^\herm$ can be viewed as a pre-processing step that can be decoupled from the NN $f_{\Theta}$ and thus does not have any impact on the learning of $\Theta$.

If at inference time \textit{only} $f_{\Theta}$ is used as the method to obtain an estimate of the ground-truth image, the method is additionally model-agnostic at inference time. The reason is that, in contrast to proper reconstruction methods that make use of the available physical model $\Au$ and its adjoint $\Au^\herm$ for reconstruction, $f_{\Theta}$ at test time ignores knowledge about the mathematical structure of the considered problem. Because of that, it is more appropriate to refer to these types of methods as post-processing or image enhancement methods, rather than reconstruction methods.

Works that focus on the reduction of artifacts and/or noise from initial images with different types of U-Nets for cardiac MRI can, for example, be found in \cite{jin2017deep} (2D U-Net), \cite{sandino2017deep} (2D U-Net with time-points stacked as channels), \cite{Hauptmann2019} (3D U-Net), or \cite{kofler2019spatio} (spatio-temporal 2D U-Net that works on $xt$- and $yt$-slices extracted from the cine MR images.

Both approaches based on the direct inversion of the forward model as well as post-processing approaches share the problem that the output of the network is in general not data-consistent. This means that these methods do not allow for the interpretation of their output as the solution of an optimization problem as for example, the one stated in \eqref{eq:mri_variational_problem_l2}.

These limitations can be overcome by embedding the aforementioned methods into proper reconstruction schemes that indeed yield solutions to minimization problems. This possibility is discussed in the next subsection.

\subsubsection{Model-Agnostic Learning of Regularization Methods}\label{sec:model_agnostic_learning_mri}

Having discussed the possibilities of NNs-based post-processing methods, one can now proceed with the integration of NNs into proper reconstruction schemes. Let us consider again the case of qualitative MR image reconstruction discussed previously, and assume we already have access to a pre-trained denoising network $f_{\Theta}$ that can reduce artifacts and noise from an initial image $\XX_0$. Given that the network $f_{\Theta}$ was trained to provide estimates of the unknown ground-truth target images for a set of MR images, we could at test time decide to utilize the output of the network, i.e.\ $\XX_{\mathrm{NN}}:=f_{\Theta}(\XX_0)$ as prior information for solving a reconstruction problem. For example, a relatively simple, yet already effective approach for 2D radial cine MRI was presented in \cite{kofler2020neural}, where for $\lambda>0$ a reconstruction problem of the type
\begin{equation}\label{eq:mri_nn_tikhonov}
    \underset{\XX}{\min}\;  \frac{1}{2}\| \Wd^{1/2} (\Ad_I \XX - \YY) \|_2^2 + \frac{\lambda}{2} \| \XX - \XX_{\mathrm{NN}}\|_2^2,
\end{equation}
is formulated and subsequently solved. Considering $\XX_{\mathrm{NN}}$ to be fixed, \eqref{eq:mri_nn_tikhonov} is strictly convex and thus has a unique solution that can be calculated by finding the image for which the gradient of the energy in \eqref{eq:mri_nn_tikhonov} with respect to $\XX$ vanishes. This coincides with solving a linear system $\Hd \XX = \mathbf{b}$ with
\begin{eqnarray}\label{eq:lin_system_for_nn_tikhonov}
\Hd &= \Au^\herm \Wd \Au +  \lambda\, \Id_N, \\
\mathbf{b} &= \Au^\herm \Wd \YY + \lambda\, \XX_{\mathrm{NN}},
\end{eqnarray}
with an appropriate iterative scheme, e.g.\ a CG method. Thus, the reconstruction  can be summarized in a two-step procedure, i.e.
\begin{eqnarray}\label{eq:nn_tikhonov_reco}
    &\XX_{\mathrm{NN}}:= &f_{\Theta}(\XX_0), \nonumber\\
    &\XX^\ast := &\underset{\XX}{\argmin}\;  \frac{1}{2}\| \Wd^{1/2} (\Ad_I \XX - \YY) \|_2^2 + \frac{\lambda}{2} \| \XX - \XX_{\mathrm{NN}}\|_2^2.
\end{eqnarray}
In \cite{kofler2020neural}, this approach was applied for radial cine MRI, where $f_{\Theta}$ was adapted from \cite{kofler2019spatio} to reduce artifacts from the initial image $\XX_0:=\Au^\herm \Wd \YY$ obtained from the density compensated $k$-space data. The work in \cite{wang2016accelerating} uses a similar approach for 2D brain MRI.

Solving problem \eqref{eq:mri_nn_tikhonov} increases data-consistency of the solution in the sense that the solution is on the one hand close to the measured $k$-space data and the other hand close to the obtained NN-prior image $\XX_{\mathrm{NN}}$. Thereby, the regularization parameter $\lambda$ controls the balance of the data-discrepancy and the proximity to $\XX_{\mathrm{NN}}$.

Depending on the exact form of the considered forward operator $\Au$, it is also possible to derive closed-form solutions for problems of the form \eqref{eq:mri_nn_tikhonov} involving a quadratic regularizer. For example, if $\Au$ denotes a single-coil Fourier operator  that samples the $k$-space data at a uniform grid, problem \eqref{eq:mri_nn_tikhonov} has an easily interpretable solution that is obtained by replacing the acquired $k$-space by a linear combination of the acquired and the one estimated from the prior $\XX_{\mathrm{NN}}$, while non-acquired $k$-space data is entirely estimated from $\XX_{\mathrm{NN}}$, see e.g.\ \cite{schlemper2017deep}.

Additionally, it is also possible to enforce strict data-consistency by solely estimating the non-acquired $k$-space coefficients from $\XX_{\mathrm{NN}}$ and keep the sampled ones unchanged, which allows the interpretation of the network as the mapping that learns the null-space of the forward model \cite{schwab2019deep}, see e.g.\ \cite{hyun2018deep} for an application to brain MRI.

From the steps in  \eqref{eq:nn_tikhonov_reco}, we see that the network $f_{\Theta}$ is only applied once to the initial image. This also means that the solution of \eqref{eq:mri_nn_tikhonov} strongly depends on the quality of the prior $\XX_{\mathrm{NN}}$. If for some reason, $\XX_{\mathrm{NN}}$ still contains artifacts and/or noise, so will the final reconstruction. 

\begin{figure}[t!]
    \centering
    \includegraphics[width=\linewidth]{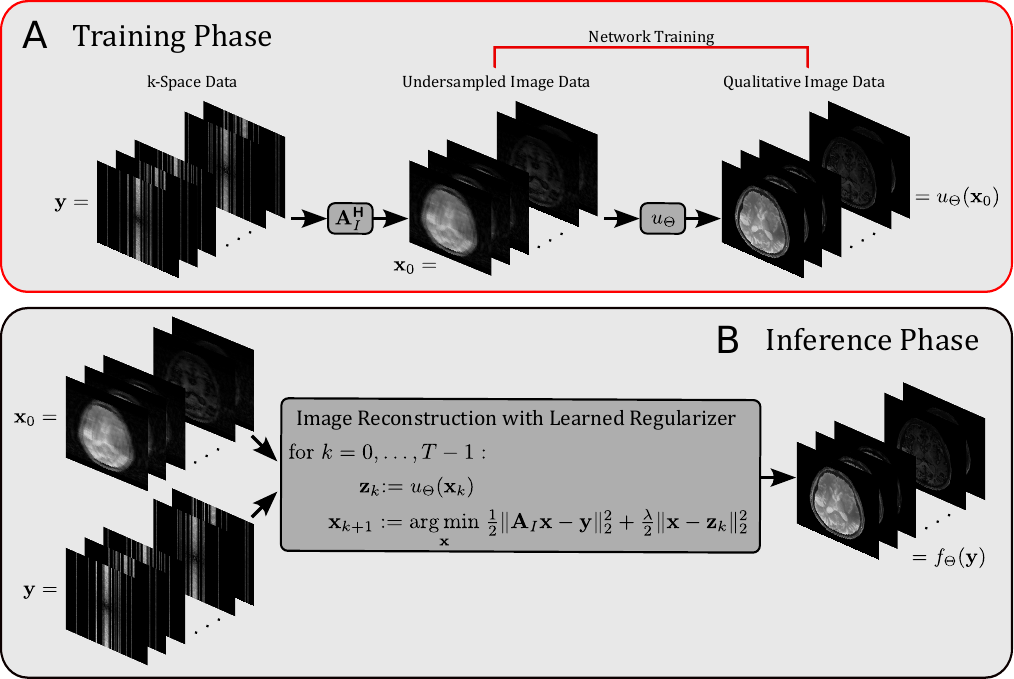}
    \caption{A schematic illustration of an example of MR image reconstruction with a model-agnostic learned regularization method based on an artifacts-removal network $u_{\Theta}$. First, in an initial step (A), the network $u_{\Theta}$ is trained to remove artifacts and noise from the initial zero-filled reconstruction. Then, the network can be used to regularize the reconstruction problem as in \eqref{eq:red_problem}. The entire mapping is here denoted by $f_{\Theta}$. Note however, that the learning of $u_{\Theta}$ is model-agnostic as \textit{first}, $u_{\Theta}$ is trained and \textit{then}, $f_{\Theta}$ is used to reconstruct the $k$-space data. }
    \label{fig:model_unaware_recon}
\end{figure}

One could heuristically opt for a repeated application of the just described two steps in \eqref{eq:nn_tikhonov_reco} to construct a sequence of solutions. However, whether or not is to possible to apply the network $f_{\Theta}$ to different images that potentially differ from the one contained in the training set 
depends on properties of $f_{\Theta}$. For example, ideally, ones has that $f_{\Theta}(\XX_{\mathrm{true}}) = \XX_{\mathrm{true}}$, i.e.\ $\XX_{\mathrm{true}}$ is a fix-point of $f_{\Theta}$, which is in general not the case.
Figure \ref{fig:model_unaware_recon} shows an example of the aforementioned model-aware reconstruction method for which the regularization method $f_{\Theta}$ is learned in a model-agnostic manner by learning to reduce noise and undersampling artifacts from an initial reconstruction.

Additionally, setting up a realistic and diverse dataset $\mathcal{D}$ that is representative enough for successful training of $f_{\Theta}$  can be challenging as it requires detailed knowledge about the specific inverse problem under consideration. Furthermore, the artifacts that the network $f_{\Theta}$ has to be able to reduce can strongly differ depending on the exact considered forward model $\Au$, making it necessary to adapt $f_{\Theta}$, for example to different acceleration factors or sampling trajectories.

This raises the question whether one can make use of general-purpose mappings $f_{\Theta}$ whose training does not rely on the availability of input-target pairs. In the signal processing community, for example, there seems to be a relatively large consensus that image denoising algorithms have started to hit the ceiling in terms of denoising performance \cite{levin2011natural}. This observation implies that one can try to take advantage from these developments and integrate image denoising algorithms as general regularizers for inverse problems.

Two methods that seek to do the above are the well-studied regularization by denoising (RED)  \cite{romano2017little} and the Plug and Play  (PnP) approach \cite{venkatakrishnan2013plug}.  PnP relies on the 
existence of a proximal map in the reconstruction algorithm (e.g. the alternating direction method of multipliers; ADMM  \cite{boyd2011distributed}) which is interpreted as a denoising step and   replaced  by a learned denoiser $u_{\Theta}$. On the other hand, RED is more general and allows to construct an explicit regularizer using a denoising network $u_{\Theta}$ to be used  with arbitrary reconstruction schemes. 

In RED, one considers a regularizing functional of the form 
\begin{equation}\label{eq:red_reg}
	\mathcal{R}(\XX) = \frac{1}{2}\,\XX^\trans \big(\XX - u_{\Theta}(\XX)\big),
\end{equation}
where  $f_{\Theta}$ is assumed to be a general Gaussian denoiser. By doing so, the regularization term in \eqref{eq:red_reg} is small either when the residual of the denoised $\XX$ is small (i.e.\ when $u_{\Theta}(\XX)\approx \XX$) or when the inner product between $\XX$ and the residual is small, which implies that the residual resembles Gaussian noise.
Under some technical assumptions on $u_{\Theta}$ (see \cite{romano2017little} and \cite{reehorst2018regularization} for details), one can show that the gradient of \eqref{eq:red_reg} can be obtained by 
\begin{equation}\label{eq:red_gradient}
\nabla_{\XX} \, \mathcal{R}(\XX) = \XX - u_{\Theta}(\XX).
\end{equation}
This means that assuming to have any pre-trained Gaussian denoiser $u_{\Theta}$ at hand, the regularizer can be employed for solving arbitrary inverse problems with gradient descent methods. In addition, $u_{\Theta}$ does not necessitate being trained on input-target pairs containing artifacts that are characteristic of the considered imaging problem, which makes the approach easily applicable to a large variety of inverse problems.
Additionally, RED was empirically shown to be applicable even for the case where the denoising method $u_{\Theta}$ does not strictly fulfill the required RED-conditions, such as state-of-the-art denoising methods based on deep CNNs, e.g.\ DnCNN \cite{zhang2017beyond}. 

A possible problem formulation for qualitative MRI can thus for example be 
\begin{equation}\label{eq:red_problem}
	\underset{\XX}{\min}\;  \frac{1}{2}\| \Ad_I \XX - \YY\|_2^2 + \frac{\lambda}{2}\, \XX^\herm \big( \XX - u_{\Theta}(\XX) \big),
\end{equation}
for which one can compute the derivative of \eqref{eq:red_problem} utilizing \eqref{eq:red_gradient} and, pursuing a fixed-point strategy, obtaining the iteration
\begin{equation}\label{eq:red_fixed_point_iters}
	\XX_{k+1} = \big(\Au^\herm  \Ad + \lambda\, \Id_N\big)^{-1}\big(\Au^\herm  \YY + \lambda\, u_{\Theta}(\XX_{k}) \big).
\end{equation}
This corresponds to a natural procedure of repeating the two steps described before, i.e.\ the application of an NN and the improvement of its output in terms of data-consistency by the solution of a minimization problem. Under the additional assumption of strong passivity \cite{romano2017little}
of $u_{\Theta}$, i.e.\ the spectral radius  of $\nabla_{\XX} u_{\Theta}$ is less or equal to one, the Hessian of $\mathcal{R}(\XX):=\XX^\herm(\XX-u_{\Theta}(\XX)$ is additionally positive definite and thus, the sequence defined by \eqref{eq:red_fixed_point_iters} converges to the solution of \eqref{eq:red_problem}.

We note that in fact, the iterative procedure in \eqref{eq:red_fixed_point_iters} corresponds to the one used to obtain the well-known and widely applied model-based deep learning (MoDL) method \cite{aggarwal2018modl} or for the deep cascade of neural networks \cite{schlemper2017deep}, where for the latter, also the NN changes across the iterations, i.e.\ one has $f_{\Theta_k}$ instead of $f_{\Theta}$. These two methods can be seen as model-aware learned versions of RED and will be discussed in more detail in the next subsection.
In addition, we point out that the aforementioned RED framework can also be adapted to incorporate artifact-removal NNs as the ones described before in Subsection \ref{sec:post_processing} instead of Gaussian denoisers, see e.g.\ \cite{liu2020rare}.

\subsubsection{Model-Aware Learning of Regularization Methods}\label{sec:model_aware}

The learning-based methods discussed in Section \ref{sec:model_agnostic_learning_mri} are reported to deliver excellent results. Additionally, they are computationally easily tractable as the learning procedure is entirely decoupled from the forward model. This is for example clearly  visible in Figures \ref{fig:full_inversion} and \ref{fig:post_processing}, where one sees that $f_{\Theta}$ never utilizes knowledge about the forward model.

For example, by comparing NNs learned in a model-agnostic manner to \eqref{eq:dico_learning_problem_L0} or \eqref{eq:cd_learning_problem_L1}, we see that until now, we have in principle only replaced the linear mappings given by the  dictionary in \eqref{eq:dico_learning_problem_L0} or the convolutional dictionary in \eqref{eq:cd_learning_problem_L1} by a more sophisticated non-linear mapping $f_{\theta}$ that is given by a non-linear NN.
However, the learning paradigm so far remained unchanged, i.e.\ we \textit{first} learned the regularizing mapping and only \textit{then} reconstructed images with an appropriate scheme.

The main focus of this section is to discuss model-aware regularization methods. While giving a precise definition of a model-aware method is rather difficult, it is fairly easy to identify a few aspects that are key and unique features of the latter, see also the remark after Table \ref{tbl:four_categories}.

Clearly, model-aware NNs only exist in the context when we are trying to predict certain quantities from observations that are generated according to a (at least partially) known physical model, e.g.\ the measured $k$-space data as in \eqref{eq:mri_forward_model} or \eqref{eq:qmri_forward_model}.
This for example explicitly excludes applications like semantic image segmentation or image classification.

The key idea of model-aware NNs is to design an NN that corresponds to a proper reconstruction scheme for an underlying image reconstruction problem. This is often achieved by so-called algorithm unrolling \cite{monga2021algorithm}, which was first introduced in the seminal work \cite{gregor2010learning}, however with a different focus. In \cite{gregor2010learning}, the aim was to accelerate a sparse coding algorithm, while these types of approaches are nowadays employed to learn entire model-aware regularization methods rather than accelerating iterative schemes. One way to derive and train model-aware reconstruction networks can be summarized as follows:
\begin{enumerate}
    \item formulate a variational reconstruction problem of the form \eqref{eq:mri_variational_problem} or \eqref{eq:qmri_variational_problem} assuming a regularization method;
    \item consider an appropriate reconstruction algorithm to solve the formulated problem;
    \item replace some components by suitable data-driven blocks that depend on learnable parameters, e.g.\ convolutional layers, deep NNs etc, and collect them all in a set of learnable parameters $\Theta$;
    \item Fix a number of iterations $T>0$, construct a network that corresponds to $T$ iterations of the derived reconstruction scheme and train according to a loss-function, e.g.\ \eqref{eq:loss_fct}.
\end{enumerate}

For example, to derive the well-known end-to-end VarNet for qualitative MR image reconstruction \cite{sriram2020end}, the starting point is the reconstruction problem
\begin{equation}\label{eq:e2e_var_net_problem}
    \underset{\XX}{\min}\;  \frac{1}{2}\| \Ad_I \XX - \YY\|_2^2 + \mathcal{R}_{\Theta}(\XX),
\end{equation}
where $\mathcal{R}_{\Theta}$ denotes an a-priori unknown learnable regularization method $\mathcal{R}_{\Theta}$. Then, to minimize problem \eqref{eq:e2e_var_net_problem} one can  use gradient descent iterations 
\begin{equation}\label{eq:e2e_var_net_scheme}
    \XX_{k+1} = g_{\YY}(\XX_k):=\XX_k - \tau \Au^\herm (\Au \XX_k - \YY) - \tau \nabla_{\XX} \mathcal{R}_{\Theta}(\XX_k).
\end{equation}
The idea is now to identify the gradient of the regularizer by a highly expressive NN, i.e.\ \ $\nabla_{\XX} \mathcal{R}_{\Theta} \equiv u_{\Theta}$, where $u_{\Theta}$, for example, denotes a U-Net. Then, we fix a number of iterations $T>0$ and  construct a model-aware NN $f_{\Theta}$ that maps the multi-coil $k$-space data to a complex-valued image by
\begin{eqnarray}\label{eq:e2e_var_net}
    f_{\Theta}:\C^{N_{\mathrm{k}} N_\mathrm{c}} &\rightarrow &\C^N \nonumber \\
               \YY &\longmapsto &\big[\big(\bigcirc_{k=1}^T g_{\YY}(\XX_k)\big) \circ \Au^\ddagger\big] (\YY),
\end{eqnarray}
where $\Au^\ddagger \in \{ \Au^\herm ,\Au^\dagger\}$ denotes some reconstruction module to obtain an initial estimate of the solution from the $k$-space data. From \eqref{eq:e2e_var_net}, we see that the network contains both learnable components (the U-Net $u_{\Theta}$) as well as model-based blocks (the gradient of the data-fidelity term $\Au^\herm (\Au \, \cdot\,  - \YY)$). Thus, the end-to-end VarNet is a model-aware NN both at training and inference time according to Table \ref{tbl:four_categories}.

As another example, VarNet \cite{hammernik2018learning} also employs a Landweber iteration to unroll the network, but uses the fields of experts model \cite{roth2005fields}
\begin{equation}\label{eq:var_net_problem}
    \mathcal{R}_{\Theta}(\XX):= \sum_{i=1}^K \langle \psi_i \big(\phi_i (\XX) \big), \mathbf{1}  \rangle,
\end{equation}
as a starting point, where $\{\psi_i\}_{i=1}^K$ and $\{\phi_i\}_{i=1}^K$ denote learnable potential functions and learnable convolutional filters.

A further widely used method, MoDL \cite{aggarwal2018modl},  can for example be obtained by identifying the proximal operator within a variable splitting method with a deep CNN. The starting point is 
\begin{equation}\label{eq:deep_cascade_problem}
    \underset{\XX,\ZZ}{\min}\;  \frac{1}{2}\| \Ad_I \XX - \YY\|_2^2 + \mathcal{R}_{\Theta}(\ZZ) + \frac{\lambda}{2} \|\XX - \ZZ\|_2^2,
\end{equation}
where the regularization is imposed on the auxiliary variable $\ZZ$ instead of $\XX$ and their equality is relaxed by the inclusion of a weighted quadratic penalty term. Then, to solve \eqref{eq:deep_cascade_problem}, one uses alternating minimization and obtains the scheme
\begin{eqnarray}\label{eq:deep_cascade_scheme}
    \ZZ_k:&= &\underset{\ZZ}{\argmin}\; \mathcal{R}_{\Theta}(\ZZ) + \frac{\lambda}{2} \|\XX - \ZZ\|_2^2 := \mathrm{prox}_{\mathcal{R}_{\Theta}}(\XX), \\
    \XX_{k+1}:&= &\underset{\XX}{\argmin}\;  \frac{1}{2}\| \Ad_I \XX - \YY\|_2^2 + \frac{\lambda}{2} \|\XX - \ZZ_k\|_2^2.
\end{eqnarray}
By setting $\mathrm{prox}_{\mathcal{R}_{\Theta}}\equiv u_{\Theta}$, this scheme allows the interpretation of the U-Net $u_{\Theta}$ as a learned denoising method. Note again the connection between the iterations in \eqref{eq:deep_cascade_problem} and the scheme derived from RED in \eqref{eq:red_fixed_point_iters}.

In principle, the possibilities to derive such reconstruction networks are endless and the number of publications that use the same underlying idea is enormous.
Further works include model-aware networks that are derived from the primal-dual hybrid-gradient (PDHG) method \cite{chambolle2011first} (originally presented for computed tomography \cite{adler2018learned} and later adopted to MRI \cite{cheng2019model}), ADMM \cite{sun2016deep, yang2018admm, yiasemis2023deep}, proximal gradient descent methods \cite{mardani2018neural}, the Neumann series \cite{gilton2019neumann} and many others, see also \cite{monga2021algorithm, shlezinger2023model} for further systematic reviews.

Algorithm unrolling-based methods have two major limitations. The first is of computational nature, as typically, the larger the number of iterations $T$ used for unrolling is, the larger the hardware requirements for training. Additionally, the number of iterations of the iterative scheme typically cannot  be varied at inference time, because convergence of the sequence generated by methods as the one in \eqref{eq:e2e_var_net_scheme} or \eqref{eq:deep_cascade_scheme}  is not guaranteed.

A family of approaches that overcomes this limitation is the one of deep equilibrium models \cite{bai2019deep}. Instead of relying on algorithm unrolling, they use implicit differentiation to derive closed-form expressions for computing the derivative of the solution mappings of fix-point iteration schemes with respect to the learnable parameters. By doing so, one can avoid the memory-intensive task of backpropagating through many iterations. For more details see \cite{bai2019deep} and \cite{gilton2021deep} for an MR application among others.

\begin{figure}[t!]
    \centering
    \includegraphics[width=\linewidth]{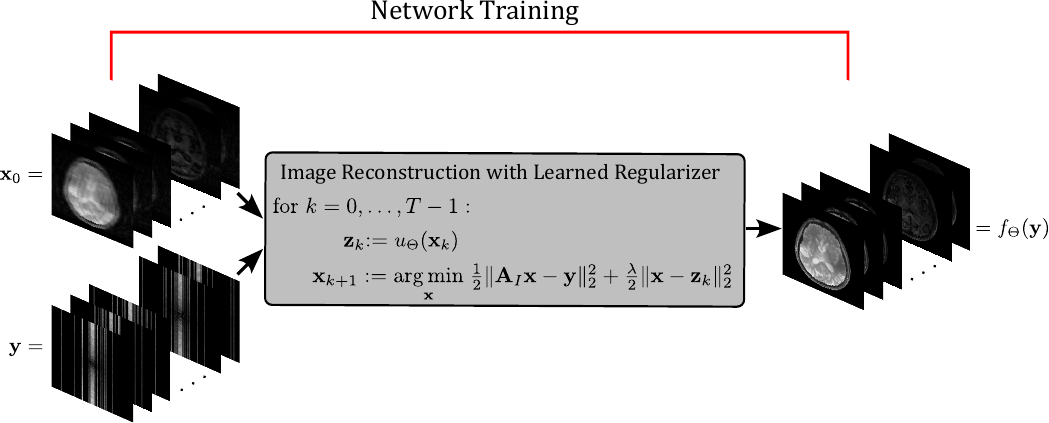}
    \caption{A schematic illustration of an example of model-aware  MR image reconstruction with a model-aware learned regularization method based on algorithm unrolling. The scheme shown in the grey box corresponds to the well-known model-based deep learning (MoDL) approach \cite{aggarwal2018modl} used for brain MR imaging or to the work used in \cite{kofler2021end} for radial cardiac cine MRI.}
    \label{fig:model_aware_recon}
\end{figure}

\begin{figure}[t!]
    \centering
    \includegraphics[width=\linewidth]{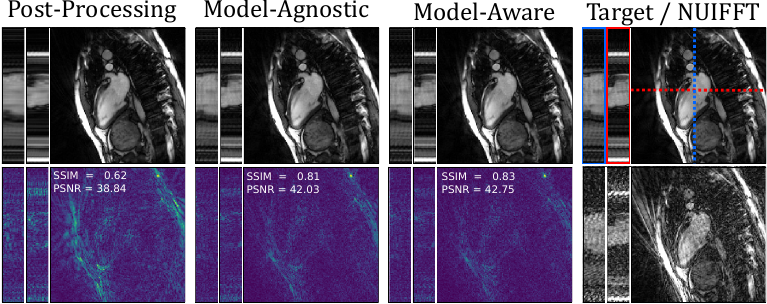}
    \caption{A comparison of a post-processing/image enhancement method, a model-unaware learned regularization method and a model-aware regularization method for accelerated dynamic cardiac MRI using radial trajecotires at an acceleration factor of approximately $9$. From left to right: CNN-based post-processing method based on the CNN-architecture resented in \cite{kofler2021end}, the CNN-basd Tikhonov regularization approach described in \eqref{eq:nn_tikhonov_reco} with the same CNN-architecture, and an end-to-end trainable version of the latter with the same CNN-architecture. The most right colum shows the target image (obtained by $kt$-SENSE \cite{tsao2003k} from which the undersampled $k$-space data was retrospectively simulated (top) and the initial reconstruction obtained by a non-uniform inverse Fourier transform (NUIFFT) applied to the density-compensated $k$-space data. Performing a proper model-based reconstruction as in \eqref{eq:nn_tikhonov_reco} improves the image compared to solely applying the post-processing CNN. Further, training the entire reconstruction method to obtain a model-aware regularization method further reduces the point-wise image-error.}
    \label{fig:post_vs_unaware_vs_aware}
\end{figure}

Figure \ref{fig:post_vs_unaware_vs_aware} shows a comparison of three different NNs-based approaches developed for radial cine MRI. The figure shows a comparison between a post-processing method, a model-agnostic reconstruction as in \eqref{eq:nn_tikhonov_reco} and it model-aware counterpart. First, to obtain the post-processing method, the CNN-architecture presented in \cite{kofler2021end} was pre-trained to reduce artifacts and noise from the initial reconstruction obtained by the applying the adjoint to the density-compensated $k$-space data. In a second step, a suitable regularization parameter $\lambda>0$ to be used in \eqref{eq:nn_tikhonov_reco} was determined on a validation dataset. Finally, the two steps in \eqref{eq:nn_tikhonov_reco} were combined into one network architecture to yield a model-aware reconstruction network that can be considered to be special case of MoDL by setting $T=1$.

\subsection{Neural Networks for  Quantitative MRI}

In contrast to qualitative MR image reconstruction, for quantitative MR reconstruction, post-processing methods that reduce noise and artifacts from an initial solution of \eqref{eq:qmri_variational_problem_splitted1} or methods to learn a regularization term for \eqref{eq:qmri_variational_problem_l2} in a model-agnostic manner  are not commonly used. Instead, most of the works either employ direct-inversion like methods that map qualitative images (which can be reconstructed with any of the previously discussed method) to estimates of the quantitative parameters, e.g.\ \cite{jeelani2020myocardial, correct, deepcest, drone, hydra} or construct model-aware reconstruction schemes, for example by unrolling a solution algorithm  for \eqref{eq:qmri_variational_problem_l2}, e.g.\ \cite{zimmermann2023pinqi, correct}.

\subsubsection{Direct Model-Inversion}\label{sec:direct_model_inversion_qmri}
In contrast to qualitative MR image reconstruction, for quantitative MR, many different methods have been proposed to obtain the quantitative parameter maps from reconstructed qualitative images by learning an approximate inverse of the signal model $q_\mathcal{M}$. The reason is that assuming the operator $\Au$ to be linear is often sufficiently accurate and thus, problems of the type \eqref{eq:mri_variational_problem_l2} are often convex, meaning that one can rely on algorithms with convergence guarantees to obtain the solution of the formulated problems. In contrast, the problem \eqref{eq:qmri_variational_problem_splitted1} can be non-convex because of the non-linearity of $q_{\mathcal{M}}$. Thus, instead of relying on gradient-descent algorithms to obtain one of possibly many solutions of \eqref{eq:qmri_variational_problem_splitted1}, one might as well employ an NN to try to invert the signal model $q_\mathcal{M}$.

The methods used for the direct inversion of $q_{\mathcal{M}}$ mainly differ in their construction and their input data. For example, 
these are either pixel-wise networks, such as MyoMapNet \cite{guo2022accelerated} which is trained to predict cardiac T1-values with fewer input images than used for the creation of the training labels, DeepCEST \cite{deepcest} (for chemical exchange saturation transfer imaging) DRONE \cite{drone}, or HYDRA \cite{hydra}. The latter two are both proposed for MRF and mainly differ in the inputs provided for the network, as HYDRA includes a (non-learned) regularized reconstruction of the qualitative images and the specifics of the 1D network. 
Further differences arise in the training, as the methods can be trained supervised with ground truth data obtained by a different method  (MyoMapNet, DeepCest) or simulated data (HYDRA, DRONE). Alternatively, as in MANTIS and the related RELAX \cite{relax}, this can either be augmented or replaced by a self-supervised loss. In both, the physical model is used only during training as part of the loss-function, not during inference.

An advantage of pixel-wise networks is that obtaining a representative sample of the data distribution is considerably easier. For training with simulated data, the real spatial correlation between ground truth labels must not be matched in the data generation, and for training with acquired data, a single acquisition leads to many input-target pairs that can be considered to be independent.
Alternatively, the network can be trained on 2D parameter maps, which allows the network to learn to make use of the spatial correlations. For example, MANTIS \cite{liu2019mantis} (presented for T2-mapping) and RCA-U-Net \cite{rcaunet} (MRF) both use U-Net-derived CNNs. These methods require either more training data, carefully chosen data augmentations and training schemes \cite{zimmermann2023ismrm} or a combination thereof.

\begin{figure}[t]
    \centering
    \includegraphics[width=\linewidth]{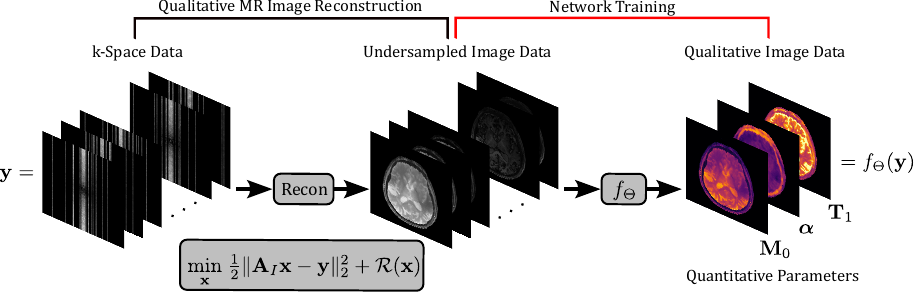}
    \caption{An example for a learned qMRI reconstruction pipeline for a T1-mapping problem. First, qualitative images have to be reconstructed with any method, learned or not learned. In a subsequent step, a network $f_{\Theta}$ is applied to the qualitative images to invert the non-linear signal model $q_{\mathcal{T}}$. Thereby, there are different choices for the network. Examples include the U-Net-like mappings \cite{liu2019mantis, rcaunet} or pixel-wise mappings \cite{guo2022accelerated, deepcest, drone, hydra}. Note that the step of qualitative image reconstruction and the training of the quantitative parameter mapping network are separated for this task. Thus, the learned $f_{\Theta}$ is neither model-aware with respect to $\Au$ nor to $q_{\mathcal{T}}$.}
    \label{fig:full_inversion_qmri}
\end{figure}

An example of a  pipeline of a T1-mapping qMRI reconstruction example is shown in Figure \ref{fig:full_inversion_qmri}. In a first step, qualitative images are reconstructed from the undersampled $k$-space data. For that, any type of method from Subsection \ref{seq:nn_for_mri} can be employed. Subsequently, a network can be trained to try to invert the model $q_{\mathcal{T}}$ to obtain an estimate of the quantitative parameters from the previously reconstructed qualitative image series.

\subsubsection{Model-Aware Learning of Regularization Methods}

Recall that in the qMR reconstruction problem \eqref{eq:qmri_forward_model} the considered operator is a composition of the non-linear signal model $q_{\mathcal{M}}$ and the data-acquisition operator $\Au$. Further, recall that an often employed strategy is to split the entire problem in two sub-problems as in \eqref{eq:qmri_variational_problem_splitted1}--\eqref{eq:qmri_variational_problem_splitted2}. This means that
different reconstruction schemes can be constructed that are model-aware only with respect to $\Au$, only with respect to $q_{\mathcal{M}}$ (rarely used) or model-aware with respect to both $\Au$ and  $q_{\mathcal{M}}$. 

A large variety of methods solve the qMR reconstruction problem as suggested in \eqref{eq:qmri_variational_problem_splitted1}--\eqref{eq:qmri_variational_problem_splitted2}, i.e.\ by first solving \eqref{eq:qmri_variational_problem_splitted2}, for example with one of the many learning-based methods described in Subsection \ref{sec:model_aware}, and then approximate a solution of \eqref{eq:qmri_variational_problem_splitted1},  with a learned direct-inversion method as described in \ref{sec:direct_model_inversion_qmri}.

For instance, DeepT1 \cite{jeelani2020myocardial} (cardiac T1-mapping) and CoRRect \cite{correct} (motion corrected $R_2^*$ mapping) both consist of two submodules, one for qualitative image reconstruction and one for parameter estimation. DeepT1 uses the model-aware reconstruction network proposed in \cite{qin2018convolutional} for the image reconstruction block, while CoRRect uses a simpler block of convolutional layers. For both methods, the acquisition model $\Au$ is used only in the qualitative image reconstruction block, whereas the parameter maps are obtained by a direct-inversion-like subnetwork based on the U-Net. For DeepT1, for example, this subnetwork is trained by using a supervised loss-function as in \eqref{eq:loss_fct_self_sup_qmri1}, where the labels were obtained from fully sampled images by a non-linear regression of the parameter model, while CoRRect employs a loss-function ot the type \eqref{eq:loss_fct_self_sup_qmri2}, i.e.\ it is trained without access to target quantitative parameters. Because of their network structure, these two methods are model-aware with respect to $\Au$ both at training as well as at inference time. DeepT1 is neither model-aware with respect to $q_{\mathcal{M}}$ at training time nor at inference time. In contrast, CoRRect is model-aware with respect $q_{\mathcal{M}}$ during training, as it employs $q_{\mathcal{M}}$ for computing the loss, but not at inference time, as the quantification of the parameters is solely based on the learned inversion without further utilizing $q_{\mathcal{M}}$. 

PGD-Net \cite{pgdnet} (MRF) unrolls a proximal gradient descent scheme on the qualitative images and includes a pre-trained, neural network-based surrogate for $q_{\mathcal{M}}$ within a learned proximal operator. Similarly, DAINTY \cite{li2021deep} uses a composition of the true $q_{\mathcal{M}}$ and a neural network as a regularizer and combines it with further low-rank and sparsity regularization of the qualitative images. Thus, both methods include the $q_{\mathcal{M}}$ as part of the regularizer for obtaining the intermediate qualitative images.  

In contrast, DOPAMINE \cite{jun2021deep} 
(also T1-mapping) includes both $\Au$ and $q_{\mathcal{M}}$ inside the data-consistency of an neural network regularized optimization problem \eqref{eq:qmri_variational_problem_l2}, which is then solved by unrolling gradient descent. Here, two subnetworks are used. One subnetwork is trained to predict from undersampled qualitative images the initial values of the parameters used in the gradient descent scheme. A second subnetwork is used to parameterize the gradient of the regularizer $\nabla \mathcal{R}$ in \eqref{eq:qmri_variational_problem_l2}. PINQI \cite{zimmermann2023pinqi} solves \eqref{eq:qmri_variational_problem_l2}  by alternating between solving an NN-regularized image reconstruction sub-problem and a parameter estimation sub-problem, both regularized by neural networks, using differentiable optimization layers \cite{bai2019deep,diffopt} that allow for end-to-end training.

\begin{figure}[t!]
    \centering
    \includegraphics[width=\linewidth]{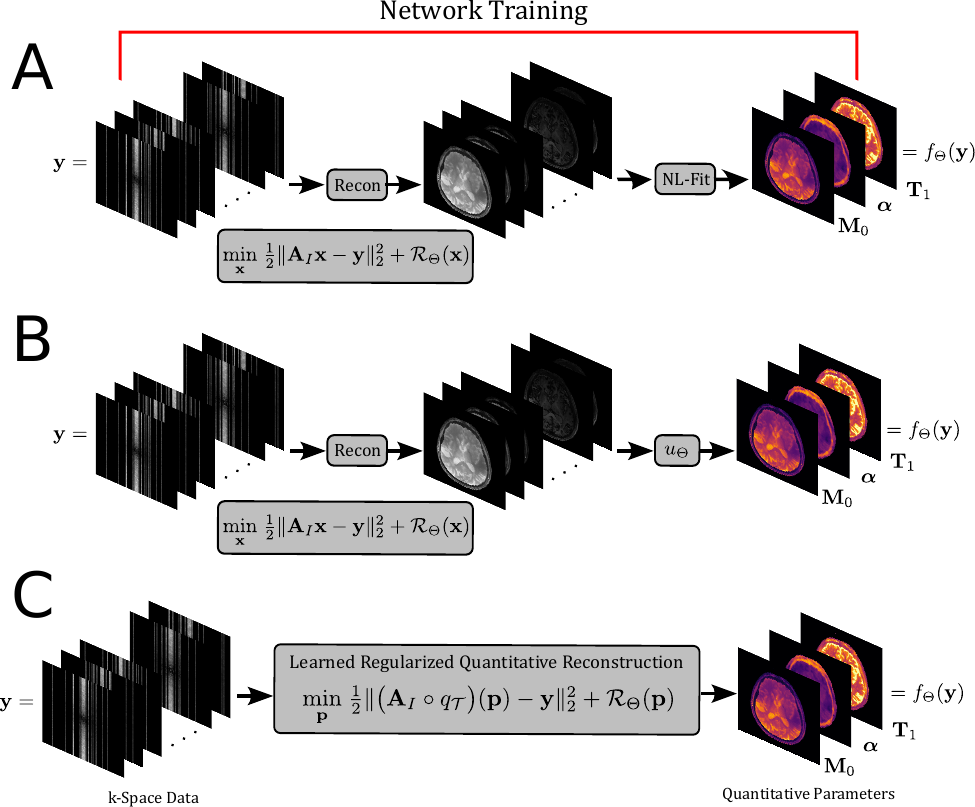}
    \caption{Three different possibile approaches for end-to-end trainable model-aware qMR reconstruction networks. In A) only the qualitative reconstruction block contains trainable parameters and a non-linear fit is utilized for the parameter quantification. In B), the non-linear fit is replaced by an entirely learned mapping $u_{\Theta}$, while C) contains a learned regularization method for the quantitative parameters.}
    \label{fig:model_aware_recon_qmri}
\end{figure}

\begin{figure}
    \centering
    \includegraphics[width=\linewidth]{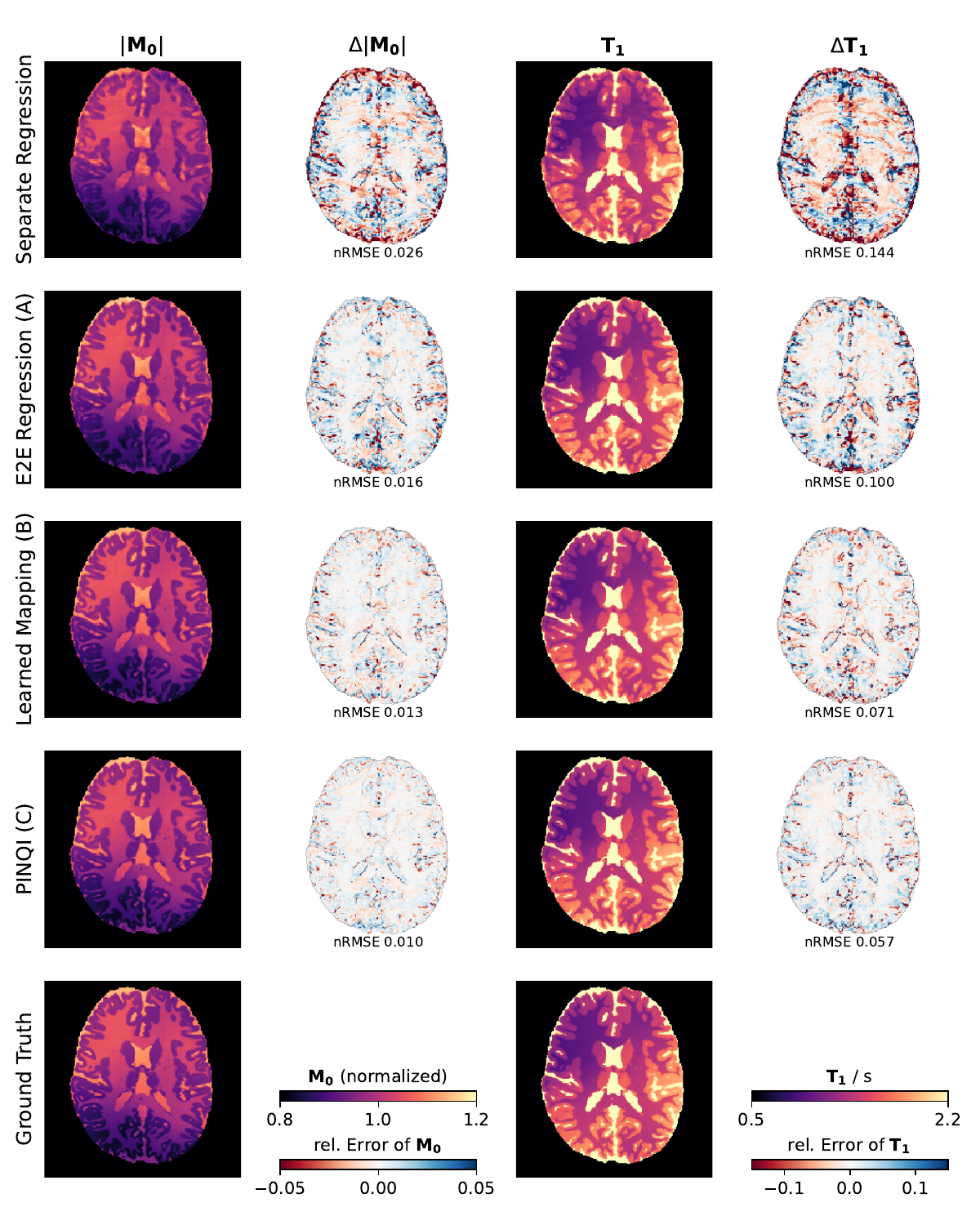}
    \caption{A comparison of examples for reconstruction of $\mathbf{M}_0$ and $\mathbf{T}_1$ from simulated 8-fold undersampled data obtained with  approaches of the types shown in Figure \ref{fig:model_aware_recon_qmri}. The method (A) of Figure \ref{fig:model_aware_recon_qmri} is shown twice. Once with a learned image reconstruction component trained in a \textit{model-agnostic} manner (top row), once trained in a  \textit{model-aware} manner (second row). Although at inference time, the procedure is exactly the same, the method trained end-to-end in a model-aware fashion outperforms the one trained in a model-agnostic way. The third row shows the results of a method corresponding to \ref{fig:model_aware_recon_qmri} B), i.e.\  a learned image reconstruction and a direct-inversion-like network for parameter estimation. 
    The fourth row shows an unrolled method for qMRI reconstruction \cite{zimmermann2023pinqi}, corresponding to \ref{fig:model_aware_recon_qmri} C), which further outperforms the first two by employing additionally learned regularization for the quantitative parameters. The last row shows the ground-truth values of $\mathbf{M}_0$ and $\mathbf{T}_1$.
    Compared to Figure \ref{fig:tv+results}, all four methods provide superior results at 8-fold acceleration.}
\label{fig:model_aware_vs_model_unaware_vs_pinqi_qmri}
\end{figure}

Figure \ref{fig:model_aware_recon_qmri} shows three different possibilities for learning end-to-end trainable model-aware qMR reconstruction methods. Figure \ref{fig:model_aware_recon_qmri} (A) resembles an end-to-end trainable version of a two-step approach as described in \eqref{eq:qmri_variational_problem_splitted1} and \eqref{eq:qmri_variational_problem_splitted2}. Thereby, only the reconstruction of the qualitative images contains a learned component, while the quantitative parameters are obtain by a standard non-linear fit, e.g.\ L-BFGS \cite{liu1989limited, nocedal1999numerical}. Note that, importantly, as indicated in the figure, the entire reconstruction method can be trained in an end-to-end fashion, for example, in a supervised manner by employing a loss-function of the type \eqref{eq:qmri_sup_training1}. This means that the regularizer $\mathcal{R}_{\Theta}(\XX)$ is trained such that the qualitative images are optimal with respect to the subsequent task, i.e.\ the non-linear fit to obtain the quantitative parameters. These characteristics make the regularization method used for the qualitative image reconstruction model-aware with respect to the operator $\Au$ and the entire reconstruction pipeline is model-aware at inference time with respect to both $\Au$ and $q_{\mathcal{T}}$.

Figure \ref{fig:model_aware_recon_qmri} (B) shows another example where the non-linear fit is replaced by a learned mapping. For example, the work in \cite{jeelani2020myocardial} falls in that category. In contrast to the method shown in Figure \ref{fig:model_aware_recon_qmri} (A), the method contains two learned components. The first one corresponds to the regularization method for the qualitative image reconstruction block, the second corresponds to a learned inversion of the model $q_{\mathcal{T}}$. Thus, both learned blocks are learned in a model-aware manner (both with respect to the model $\Au)$, but only the first reconstruction block is model-aware at inference time, since the quantification of the parameters is performed in a model-agnostic way by $u_{\Theta}$.

Finally, Figure \ref{fig:model_aware_recon_qmri} (C) shows an end-to-end trainable qMR reconstruction method that employs a learned regularization that is applied to the quantitative parameters. The work in \cite{zimmermann2023pinqi} (PINQI), presents an approach  that alternates between the reconstruction of regularized qualitative images and a regularized non-linear fit to obtain the quantitative parameters. Because of the specific construction of the network and the end-to-end trainability, PINQI represents a method that employs model-aware regularization components and is additionally also model-aware at inference time.

Note that the method in Figure \ref{fig:model_aware_recon_qmri} (A) strongly differs from its model-agnostic counterpart, where one would first train a reconstruction network to yield qualitative images, as for example shown in Figure \ref{fig:model_aware_recon}, and only in a second step perform a non-linear fit. Figure \ref{fig:model_aware_vs_model_unaware_vs_pinqi_qmri} shows a comparison of the two different just described approaches, which are parts of an ablation study in the work \cite{zimmermann2023pinqi}. Clearly, the model-aware approach (left) surpasses the model-agnostic approach (right) in terms of reconstruction accuracy for the obtained quantitative parameters.

\subsection{Scan-Specific Networks}

A further quite successful possibility to employ NNs for MR image reconstruction, even in the total absence of training data, is the use of so-called \textit{untrained} or \textit{scan-specific networks} such as the well-known Deep Image Prior (DIP) \cite{ulyanov2018deep}, the Deep Decoder (DD)\cite{heckel2018deep}, or certain Implicit Neural Representation (INR) \cite{sitzmann2020implicit} methods. In these methods, the networks are optimized for one specific input. Thus it is not possible to identify a separate training and inference phase. Instead, these approaches often use a reparametrization of the image by an NN. In DIP, for example, given one specific set of acquired data $\YY$, the energy of the following  minimization problem 
\begin{equation}\label{eq:dip_problem}
    \underset{\Theta}{\min}\;  \frac{1}{2}\| \Ad_I f_{\Theta}(\ZZ) - \YY\|_2^2,
\end{equation}
is decreased using some gradient descent scheme combined with a suitable early stopping criterion or some additional regularization, see next paragraph. Here $\ZZ$ is a fixed input that is typically chosen as random noise and $f_{\Theta}$ is an NN with parameters $\Theta$. By comparing \eqref{eq:dip_problem} to \eqref{eq:mri_variational_problem_l2}, we see that, instead of minimizing over the image, ones seeks for the network parameters that minimize the error between the measured and predicted data. Once \eqref{eq:dip_problem} the energy is decreased and the stopping criterion is satisfied, the image can be retrieved by $\XX^\ast:=f_{\Theta}(\ZZ)$.
Note that, minimizing \eqref{eq:dip_problem} conceptually coincides with the procedure of \textit{training} an NN $f_{\Theta}$ with a self-supervised loss-function as in \eqref{eq:loss_fct_self_sup_mri} with only $N_{\mathrm{train}}=1$ sample.

Note that for the case where $f_{\Theta}$ is a deep and overparametrized NN, strictly minimizing \eqref{eq:dip_problem} typically can result to overfitting of the acquired data and a noisy or artifact-distorted result. Hence, additional regularization is required, either by early stopping as suggested in the original work \cite{ulyanov2018deep}, by including further regularization terms \cite{sitzmann2020implicit} or by constraining the NN topology \cite{heckel2018deep}. 
Scan-specific methods have successfully been applied for MR image reconstruction, e.g.\ non-Cartesian dynamic cardiac MRI \cite{yoo2021time,catalan2023unsupervised}, for static 2D MRI \cite{arora2020untrained, darestani2021accelerated, feng2023scan}, and for qMRI \cite{heydari2024joint,gao2021accurate}. For image reconstruction, the approaches were reported to achieve similar results as post-processing CNN-based methods trained with target data in a supervised manner, but are typically outperformed by unrolled methods such as the End-to-end-Varnet \cite{sriram2020end}.

\subsection{Summary}
Here, we provide a short summary of some of the aforementioned methods. Table \ref{tbl:works_classification} lists some of the works and classifies them according to Table \ref{tbl:four_categories} depending on whether or not some of the components are learning-based as well as whether the learning-paradigm is model-aware or not. Note that, as we have previously seen, qMRI problems are often tackled by addressing two different problems. In that case, for qMRI problems, there are three different possibilities to employ learned regularization; either only in for the qualitative MR problem, only for the quantitative MR problem or for both.

\begin{table}[h!]
\caption{A summary of different data-driven methods for qualitative and quantitative MR image reconstruction which differ from each other with respect to the employed training paradigm. For the different approaches, a categorization between model-aware and model-agnostic learning is performed. Further, to emphasize whether the learning of the regularization method is model-aware or not, the distinction is done both for the training phase as well as for the inference phase. We should note however that this categorization  is not strictly objective as the boundaries of the different categories are not always rigid. We use the following abbreviations: Radial (R); Cartesian (C), Self-supervised (Self), Diffusion Tenson Imaging (DTI).
}
\setlength\extrarowheight{1pt}
\resizebox{\columnwidth}{!}{
\begin{tabular}{l |c c |c c |c}
    \hline
    \rowcolor{gray!30}
    \textbf{Name and reference} & \multicolumn{2}{c|}{\textbf{Qual.\ Recon.\ Block}}  &  \multicolumn{2}{c|}{\textbf{Quant.\ Recon.\ Block}} & \textbf{Notes}\\
    & \multicolumn{2}{c|}{\cellcolor{gray!20}Model-Awareness} & \multicolumn{2}{c|}{ \cellcolor{gray!20} Model-Awareness} & \\
    & Training & Inference & Training & Inference & \\  
    3D U-Net \cite{Hauptmann2019} & no & no & N.A. & N.A. & Dynamic Cardiac (R)\\
    XT,YT \cite{kofler2019spatio} & no & no & N.A. & N.A. & Dynamic Cardiac (R)\\
    CNN-Tikhonov \cite{kofler2020neural} & no & yes & N.A. & N.A. & Dynamic Cardiac  (R)\\
    Null-Space learning \cite{hyun2018deep} & yes & yes & N.A. & N.A. & Brain \\
    E2E CNN-Tikhonov \cite{kofler2021end} & yes & yes & N.A. & N.A. & Dynamic Cardiac  (R)\\
    MoDL \cite{aggarwal2018modl} & yes & yes & N.A. & N.A. &  Brain \\
    Deep Cascades \cite{aggarwal2018modl, qin2018convolutional} & yes & yes & N.A. & N.A. & Dynamic Cardiac   (C)\\
    VarNet \cite{hammernik2018learning} & yes & yes & N.A. & N.A. & Knee  (C)\\ 
    2D DLMRI \cite{ravishankar2010mr} & no & yes & N.A. & N.A. & Brain  (C)\\
    DLMRI + TV \cite{caballero2014dictionary, wang2013compressed} & no & yes & N.A. & N.A. & Dynamic Cardiac   (C)\\
    CDL \cite{quan2016compressed} & no & yes & N.A. & N.A. & Dynamic Cardiac  (C)\\
    NN-CDL \cite{kofler2022cdl} & yes & yes & N.A. & N.A. & Dynamic Cardiac   (R)\\
    NN-CAOL \cite{kofler2022caol} & yes & yes & N.A. & N.A. & Dynamic Cardiac   (R)\\
    Deep DLMRI \cite{kofler2022nndlmri} & yes & yes & N.A. & N.A. & Dynamic Cardiac   (C)\\ 
    qDLMRI \cite{kofler2023quantitative} & N.A. & yes & no & yes & T1-mapping, Brain  (R)\\ 
    NoSENSE (qMRI) \cite{zimmermann2023nosense} & yes & yes & N.A. & yes & T1-mapping, Cardiac  (C)\\
    DDM2 \cite{ddm2} & no & yes & N.A. & yes & DTI\\
    DeepT1\cite{jeelani2020myocardial} & yes & yes & no & no & T1-mapping, Cardiac (C)\\
    PnP for QMRI \cite{fatania2022plug} & no & yes & N.A. & yes & MRF, Brain\\
    DeepCest \cite{deepcest} &  N.A. & N.A. & no & no & Inversion of $q_{\mathcal{M}}$\\ 
    MyoMapNet \cite{guo2022accelerated} &  N.A. & N.A. & no & no & Inversion of $q_{\mathcal{M}}$\\ 
    DRONE \cite{drone} & N.A. & N.A. & no & no & MRF (pixelwise)\\ 
    RCA-U-Net \cite{rcaunet} & N.A. & N.A. & no & no & MRF (CNN)\\ 
    MANTIS \cite{liu2019mantis} & yes & N.A. & yes & no & supervised and self\\ 
    RELAX \cite{relax} & yes & N.A. & yes & no & T2-mapping, (self)\\
    PGD-Net \cite{pgdnet} & yes & yes &  yes & no & Only surrogate of $q_{\mathcal{M}}$ \\ 
    DAINTY \cite{li2021deep} & yes & yes &  yes & yes & T1-mapping, Brain  (C)\\   
    CoRRECT \cite{correct} & yes & yes & yes & no & R$_2^*$-mapping, Brain\\ 
    DOPAMINE \cite{jun2021deep} & yes & yes & yes & yes & Unrolled gradient descent\\ 
    PINQI \cite{zimmermann2023pinqi} & yes & yes & yes & yes & Alternating optimzation\\ 
\end{tabular}
}
\label{tbl:works_classification}
\end{table}

Note that the Table is not exhaustive concerning all published methods. Instead, it may serve as guidance to locate relevant papers, providing a starting point for delving deeper into the topic.

\section{Conclusions and Outlook}

We conclude this chapter with a brief discussion about important aspects that must be taken into consideration when applying NNs for qualitative and quantitative MR image reconstruction.

The most important limitation of NNs, despite their empirical excellent performance, is the fact that from a mathematical point of view, they are black-boxes and they lack many of the advantages of the model-based methods listed at the end of Subsection \ref{sec:reg_qualitative}. Some methods aim to give some sort of explanation of their functioning mechanisms, e.g.\ in the form of relevance maps that provide insight into the area of the input image that determines the NN's decision, but these are mainly restricted to classification tasks, see e.g.\ \cite{montavon2018methods}.
An alternative way to maintain  interpretability, is to strongly constrain the family of learnable regularizers. For example, the well-known VarNet \cite{hammernik2018learning} learns convolutional filters by unrolling a reconstruction algorithm that minimizes a functional for which the regularizer stems from the field-of-experts model \cite{roth2005fields}. As such, one obtains an explicit regularization method and the reconstruction scheme can be linked to an underlying functional. Other works that are also based on algorithm unrolling, e.g.\ \cite{kofler2022cdl, kofler2022caol}, further constrain the regularization method and learn unit-norm convolutional filters that can be associated with synthesis and analysis sparsifying transforms, i.e.\ the regularizing mechanism of the obtained method is entirely interpretable. Similarly, patch-based sparsity-transforms as the ones discussed in \ref{sec:data_driven_mri} can as well be learned \cite{ravishankar2017physics, kofler2022nndlmri}.

In addition, it is also possible to choose a hand-crafted regularization method, e.g.\ TV or Wavelets, and use deep NNs to enhance the performance of the latter. For example, in \cite{Kofler_SIAM_Imaging_2023}, regularization parameter maps as the ones described in \eqref{def:TV_spatial} were learned with a U-Net by unrolling the primal-dual hybrid-gradient method (PDHG) \cite{chambolle2011first}. As such, the estimation of the regularization parameter-maps is not interpretable because of the black-box of the U-Net, but the subsequent reconstruction is, because its corresponds to a convex minimization method with convergence guarantees.
Analogously, in \cite{nguyen2023map}, scale-dependent regularization parameters were learned for a wavelet-based regularization.

In general, the black-box character of deep NNs poses a challenge for their clinical application because of the existence of instabilities and so-called adversarial attacks, where small perturbations of the input can severely change the output\cite{adversial,antun2020instabilities}. In addition, issues related to the transferability of NNs-based methods to datasets that differ from the one used for training remain a challenge. For example, \cite{bhadra2021hallucinations} reported examples of realistically appearing hallucinations that the reconstruction method yielded when tested on MR images of children while being trained on MR images of solely adults.

To address these issues, there exist several different possibilities. First of all, the consensus seems to be that methods based on the direct inversion of the forward model, as presented in \cite{zhu2018image} should be preferably avoided. Including prior knowledge about the reconstruction problem by making use of the physical model, both at inference time and/or during training, seems to at least partially mitigate this issue. The results reported in \cite{antun2020instabilities} suggest that the direct inversion method AUTOMAP \cite{zhu2018image} was the least stable concerning adversarial attacks, while post-processing methods as DAGAN \cite{yang2017dagan} or unrolled methods as the deep cascade of CNNs \cite{schlemper2017deep} and VarNet \cite{hammernik2018learning}  seem to show somewhat less drastic artifacts. 

In contrast, the authors in \cite{genzel2022solving} reported that end-to-end trained reconstruction methods, if trained properly, can be robust with respect to both noise and adversarial perturbations. They emphasized the importance of avoiding inverse crimes during training, stressing the importance of techniques like jittering\cite{bishop1995training}, i.e.\ a perturbation of the input by random noise.  In particular, the authors reported that their considered networks were not more unstable than the non-learning method based on TV-minimization which was unrolled using ADMM \cite{glowinski1975approximation, gabay1976dual}.

Further, the impact of distribution-shifts was investigated in \cite{darestani2021measuring} both for sparsity-based methods employing wavelet transforms \cite{chen2012compressive}, untrained networks \cite{darestani2021accelerated} as well as unrolled methods \cite{sriram2020end}. The authors concluded that all methods comparably suffer from distribution shifts. However, for $\ell_1$-minimization based methods such as TV or wavelets, the performance gap can be easily closed by re-adjusting the only tunable parameter, i.e.\ the regularization parameter $\lambda$, to the newly considered dataset. In contrast, adapting pre-trained networks to new datasets is a non-trivial task that can be addressed with different techniques. These include model adaptation \cite{gilton2021model, alanov2022hyperdomainnet, kanakis2020reparameterizing} or test-time training \cite{darestani2022test}.

Last, we note that the training of state-of-the-art methods based on deep learning models can often take several days. In addition, there is often the need to repeat the training procedure for different choices of hyper-parameters (e.g.\ number of layers, number of filters for different layers, different learning rates, etc.). Unrolled methods with deeper networks were observed to yield superior results compared to their counterparts with more shallow networks \cite{kofler2022more}. Additionally, most of the winning teams in image reconstruction challenges nowadays seem to employ methods based on algorithm unrolling \cite{muckley2021results, beauferris2022multi, sidky2022report},
which in contrast to model-agnostic methods, are linked to substantial hardware requirements. These aspects raise questions about the large carbon footprint of deep learning-based methods in general as well as about their inclusiveness among different research institutions and countries \cite{xu2021survey}.

Despite these open questions, deep learning-based methods have clearly emerged as a powerful and transformative approach in the field of MR image reconstruction and have yielded unprecedented results. Their remarkable success in this domain potentially holds the promise of enhancing diagnostic accuracy and advancing the field of both qualitative and quantitative MR.

\bibliographystyle{plain}
\bibliography{bib_jabref_noduplicates}

\end{document}